\documentclass[11pt]{article}
\usepackage{amssymb}

\usepackage{amsmath}
\usepackage{booktabs}
\usepackage{graphicx}
\usepackage{theorem}

\newtheorem{thm}{Theorem}
\newtheorem{lem}{Lemma}

\theorembodyfont{\rmfamily}

\allowdisplaybreaks[3]
\makeatletter
\newcommand{\Proofname}{Proof}
\makeatother
\makeatletter
\def\BOXSYMBOL{\RIfM@\bgroup\else$\bgroup\aftergroup$\fi
  \vcenter{\hrule\hbox{\vrule height.85em\kern.6em\vrule}\hrule}\egroup}
\makeatother
\newcommand{\BOX}{%
  \ifmmode\else\leavevmode\unskip\penalty9999\hbox{}\nobreak\hfill\fi
  \quad\hbox{\BOXSYMBOL}}

\def\hbbe{{\widehat{\bbe}}}

\def\al{{\alpha}}

\def\de{{\delta}}
\def\ep{{\varepsilon}}
\def\la{{\lambda}}

\def\bbe{{\text{\boldmath $\beta$}}}

\def\bep{{\text{\boldmath $\varepsilon$}}}

\def\bth{{\text{\boldmath $\theta$}}}

\def\bro{{\text{\boldmath $\rho$}}}
\def\bpsi{{\text{\boldmath $\psi$}}}

\def\bel{{\text{\boldmath $\ell$}}}

\def\bbeh{{\widehat \bbe}}

\def\bthh{{\widehat \bth}}

\def\bbet{{\widetilde \bbe}}

\def\btht{{\widetilde \bth}}

\def\Si{{\Sigma}}

\def\La{{\Lambda}}

\def\bSi{{\text{\boldmath $\Si$}}}

\def\bLa{{\text{\boldmath $\La$}}}

\def\bPsi{{\text{\boldmath $\Psi$}}}

\def\bSih{{\widehat \bSi}}
\def\bPsih{{\widehat \bPsi}}

\def\bLah{{\widehat \bLa}}

\def\a{{\text{\boldmath $a$}}}

\def\c{{\text{\boldmath $c$}}}

\def\v{{\text{\boldmath $v$}}}

\def\y{{\text{\boldmath $y$}}}

\def\A{{\text{\boldmath $A$}}}
\def\B{{\text{\boldmath $B$}}}

\def\D{{\text{\boldmath $D$}}}

\def\G{{\text{\boldmath $G$}}}
\def\H{{\text{\boldmath $H$}}}
\def\I{{\text{\boldmath $I$}}}
\def\J{{\text{\boldmath $J$}}}

\def\P{{\text{\boldmath $P$}}}
\def\Q{{\text{\boldmath $Q$}}}

\def\W{{\text{\boldmath $W$}}}
\def\X{{\text{\boldmath $X$}}}

\def\yb{{\overline \y}}

\def\Xb{{\overline \X}}

\def\ybt{{\tilde \y}}

\def\vbt{{\tilde \v}}

\def\Xbt{{\widetilde \X}}

\def\Nc{{\cal N}}

\def\Wc{{\cal W}}

\def\Re{{\mathbb{R}}}
\def\tr{{\rm tr\,}}
\def\diag{{\rm diag\,}}

\def\[{{\text{\boldmath $[$}}}
\def\]{{\text{\boldmath $]$}}}
\def\et{{\it et\, al.}}
\def\zero{{\bf\text{\boldmath $0$}}}
\def\one{{\bf\text{\boldmath $1$}}}

\def\|{{\,|\,}}

\def\/{{\Bigr/\!\!}}

\def\1r{{\rm (1)}}
\def\2r{{\rm (2)}}
\def\3r{{\rm (3)}}
\def\4r{{\rm (4)}}
\def\5r{{\rm (5)}}

\def\non{{\nonumber}}
\def\one{{\bf\text{\boldmath $1$}}}
\def\zero{{\bf\text{\boldmath $0$}}}

\def\vbt{\widetilde{\bf\text{\boldmath $v$}}}
\def\ybt{\widetilde{\bf\text{\boldmath $y$}}}

\begin{document}
\title{Empirical Best Linear Unbiased Predictors in Multivariate Nested-Error Regression Models}

\author{
Tsubasa Ito\footnote{Graduate School of Economics, University of Tokyo, 7-3-1 Hongo, Bunkyo-ku, Tokyo 113-0033, JAPAN. {E-Mail: tsubasa$\_$ito.0710@gmail.com}}
and
Tatsuya Kubokawa\footnote{Faculty of Economics, University of Tokyo, 7-3-1 Hongo, Bunkyo-ku, Tokyo 113-0033, JAPAN. \newline{E-Mail: tatsuya@e.u-tokyo.ac.jp }} 
}

\date{}
\maketitle
\begin{abstract}
For analyzing unit-level multivariate data in small area estimation, we consider the multivariate nested error regression model (MNER) and provide the empirical best linear unbiased predictor (EBLUP) of a small area characteristic based on second-order unbiased and consistent estimators of the ^^ within' and ^^ between' multivariate components of variance.
The second-order approximation of the mean squared error (MSE) matrix of the EBLUP and its unbiased estimator are derived in closed forms.
The confidence interval with second-order accuracy is also provided analytically.

\par\vspace{4mm}
\end{abstract}

\section{Introduction}

Linear mixed models and model-based predictors in small area estimation have been studied extensively and actively in recent years due to the growing demand for reliable small area estimates.
In small area estimation, direct design-based estimates for small area means have large standard errors due to small sample sizes from small areas.
In order to improve accuracy, the linear mixed models are considered which consist of fixed effects based on common parametes and random effects depending on areas, and the resulting empirical best linear unbiased predictors (EBLUP) provide more reliable estimates by ^^ borrowing strength' from neighboring areas.
The linear mixed models used in small area estimation are the Fay-Herriot model for analyzing area-level data by Fay and Herriot (1978) and the nested error regression (NER) model for analyzing unit-level data by Battese, Harter and Fuller (1988).
Various extensions and generalizations of these models and many statistical methods for inference  have been studied in the literature.
For comprehensive reviews of small area estimation, see Ghosh and Rao (1994), Datta and Ghosh (2012), Pfeffermann (2013) and Rao and Molina (2015).

\medskip
In this paper, we consider multivariate nested error regression models (MNER) with fixed effests based on a vector of regression coefficients $\bbe$ and vectors of random effects $\v_i$ and sampling errors $\bep_{ij}$ for the $j$-th unit in the $i$-th area.
When $\bth_a$, defined by $\bth_a=\c_a^\top\bbe+\v_a$, is a characteristic of interest for the $a$-th area and constant $\c_a$, the Bayes estimator of $\bth_a$ in the Bayesian context is 
$$
\btht_a(\bbe, \bPsi, \bSi)=\c_a^\top\bbe + \B_a (\yb_a-\Xb_a^\top\bbe),
$$
for $\B_a=\bPsi(\bPsi+n_a^{-1}\bSi)^{-1}$, where $\bSi$ and $\bPsi$ are covariance matrices of $\bep_{ij}$ and $\v_i$, respectively, $n_a$ is a size of a sample from the $a$-th area, and $\yb_a$ and $\Xb_a$ are sample means of response variables and the associated explanatory variables in the $a$-th area.
When components of $\v_i$ and $\bep_{ij}$ are mutually independent, namely $\bPsi$ and $\bSi$ are diagonal matrices, it is enough to treat the estimation of each component of $\bth_a$ separately.
When components of $\v_i$ or $\bep_{ij}$ are correlated each other, however, it could be better to consider the estimation of $\bth_a$ simultaneously.
For example, the survey and satellite data of Battese, $\et$ (1988) consist of two crop areas under corn and soybean, and it should be reasonable that the two crop areas are correlated each other.

\medskip
The multivariate small area estimation has not been studied so much, while most results in small area estimation have been provided in the univariate cases.
Fay (1987) proposed a multivariate Fay-Herriot model for analyzing multivariate area-level data.
Porter, Wikle and Holan (2015) and Benavent and Morales (2016) suggested multivariate spatial Fay-Herriot models with covariance matrices in which spatial dependence is embedded.
Concerning the multivariate nested error regression (MNER) models, Fuller and Harter (1987) obtained the empirical Bayes estimator or the EBLUP and the analytical results for its uncertainty, and Datta, Day and Maiti (1998) developed the fully Bayesian approach.
Datta, Day and Basawa (1999) also provided general theoretical results for the multivariate empirical Bayes estimators, but did not give concrete expressions in the fully unknown case of covariance matrices.

\medskip
The MNER model has the two components of covariance: ^^ between' component $\bPsi$ and ^^ within' component $\bSi$.
We here use an exact unbiased estimator $\bSih$ for $\bSi$, and  for $\bPsi$, we suggest a nonnegative definite and consistent estimator $\bPsih$  which is a second-order unbiased estimator of $\bPsi$.
For the other estimation methods, see Calvin and Dykstra (1991a, b).
Substituting $\bPsih$ and $\bSih$ into $\bPsi$ and $\bSi$ in the Bayes estimator $\btht_a(\bbe,\bPsi, \bSi)$ and estimating $\bbe$ by the generalized least squares estimator $\bbeh$, one gets the empirical Bayes estimator or EBLUP $\bthh_a^{EB}=\btht_a(\bbeh,\bPsih,\bSih)$.
We derive analytically a second-order approximation of the MSE matrix of the EBLUP and provide a closed form expression of a second-order unbiased estimator, denoted by ${\rm msem}(\bthh_a^{EB})$, of the MSE matrix of the EBLUP.
These results are extensions of the univariate case.
It is noted that similar results were given by Fuller and Harter (1987) who considered to estimate $\B_a$ nearly unbiasedly, which is slightly different from the approach of this paper.

\medskip
Another topic addressed in the paper is the confidence interval problem.
As pointed out in Diao, Smith, Datta, Maiti and Opsomer (2014), one difficulty with traditional confidence intervals is that the coverage probabilities do not have second-order accuracy.
It is also numerically confirmed that  the coverage probabilities are smaller than the nominal confidence coefficient.
Diao et al. (2014) suggested the construction of accurate confidence interval based on the EBLUP and the estimator of MSE of EBLUP so that the coverage provability is correct up to second order.
For other studies on the confidence interval problem, see Datta, Ghosh, Smith, Lahiri (2002), Basu, Ghosh and Mukerjee (2003), Chatterjee, Lahiri and Li (2008), Kubokawa (2010),  Sugasawa and Kubokawa (2015) and Yosimori and Lahiri (2014).
In this paper, we consider the confidence interval for the liner combination $\bel^\top\bth_a$ for $\bel \in \mathbb R^k$.
The naive confidence interval is given by $\bel^\top \bthh_a^{EB}\pm z_{\al/2}\times\sqrt{\mathstrut \bel^\top{\rm msem}(\bthh_a^{EB}) \bel}$ where $z_{\al/2}$ is the $100(1-\al/2)\%$ percentile of the standard normal distribution.
Because this confidence interval does not have second-order accuracy, using similar aguments as in Diao et al. (2014), we construct the closed-form confidence interval whose coverage probability is identical to the nominal confidence coefficient $1-\al$ up to second order.

\medskip
The paper is organized as follows:
In Section \ref{sec:EBLUP}, we probide an exact unbiased estimator $\bSih$ for $\bSi$ and a nonnegative definit, consistent and second-order unbiased estimator $\bPsih$ of $\bPsi$.
Substituting these estimators into the Bayes estimator yields the empirical Bayes estimator or EBLUP $\btht_a(\bbeh,\bPsih,\bSih)$.
In Section \ref{sec:MSE}, we derive a second-order approximation of the MSE matrix of the EBLUP and a second-order unbiased estimator of the MSE matrix analytically.
Section \ref{sec:ci} presents the confidence interval with second-order accuracy..
The numerical investigation and emprical studies are given in Section \ref{sec:sim}.

\section{Empirical Best Linear Unbiased Prediction}
\label{sec:EBLUP}

In this paper, we assume that data $(\y_{ij}, \X_{ij})$ for $i=1\ldots,m$ and $j=1,\ldots,n_i$ are observed, where $m$ is the number of small areas, $n_i$ is the number of the subjects in an $i$-th area such that $\sum_{i=1}^m n_{i}=N$, $\y_{ij}$ is a $k$-variate vector of direct survey estimates and $\X_{ij}$ is a $s \times k$ matrix of covariates associated with $\y_{ij}$ for the $j$-th subject in the $i$-th area.
Then, we assume the multivariate nested-error regression model described as
\begin{equation}
\y_{ij} = \X_{ij}^\top \bbe + \v_i + \bep_{ij},
\quad i=1, \ldots, m, 
\quad j=1, \ldots, n_i,
\label{eqn:MNE}
\end{equation}
where $\bbe$ is an $s$-variate vector of unknown regression coefficients, $\v_i$ is a $k$-variate vector of random effects depending on the $i$-th area and $\bep_{ij}$ is a $k$-variate vector of sampling errors.
It is assumed that $\v_i$ and $\bep_{ij}$ are mutually independently distributed as
$$
\v_i \sim \Nc_k (\zero, \bPsi)\quad \text{and}\quad \bep_{ij}\sim\Nc_k(\zero, \bSi),
$$
where $\bPsi$ and $\bSi$ are $k\times k$ unknown and nonsingular covariance matrices.

\medskip
We now express model (\ref{eqn:MNE}) in a matrix form.
Let $\y_i=(\y_{i1}^\top, \ldots, \y_{in_i}^\top)^\top$, $\y=(\y_1^\top, \ldots, \y_m^\top)^\top$, $\X_i=(\X_{i1},\ldots,\X_{in_i})^\top$, $\X=(\X_1^\top, \ldots, \X_m^\top)^\top$, $\bep_i=(\bep_{i1}^\top, \ldots, \bep_{in_i}^\top)^\top$ and $\bep=(\bep_1^\top, \ldots, \bep_m^\top)^\top$.
Then, model (\ref{eqn:MNE}) is expressed as
\begin{equation}
\y_i=\X_i\bbe + \one_{n_i} \otimes \v_i + \bep_i,
\label{eqn:MNE1}
\end{equation}
where $\one_{n_i} \otimes \v_i \sim \Nc_{kn_i}(\zero, \J_{n_i} \otimes \bPsi)$ and $\bep_i\sim\Nc_{kn_i}(\zero, \I_{n_i} \otimes \bSi)$ for $\J_{n_i}=\one_{n_i}\one_{n_i}^\top$.

\medskip
For the $a$-th area, we want to predict the quantity $\bth_a=\c_a^\top \bbe+\v_a$, which is the conditional mean $E[\y_a\mid \v_a]$ given $\v_a$ when 
$$
\c_a=\Xb_a=n_a^{-1}\sum_{j=1}^{n_a}\X_{aj}.
$$
A reasonable estimator can be derived from the conditional expectation $E[\bth_a\mid \y_a]=\c_a^\top \bbe+E[\v_a \mid \y_a]$.
The conditional distribution of $\v_i$ given $\y_i$ and the marginal distribution of $\y_i$ are
\begin{equation}
\begin{split}
\v_i \mid \y_i \sim & \Nc_k(\vbt_i(\bbe, \bPsi, \bSi), (\bPsi^{-1}+n_i\bSi^{-1})^{-1}),\\
\y_i \sim& \Nc_{kn_i}(\X_i\bbe, \J_{n_i}\otimes\bPsi+\I_{n_i}\otimes\bSi),
\end{split}
\quad i=1, \ldots, m,
\label{eqn:post}
\end{equation}
where
\begin{equation}
\vbt_i(\bbe,\bPsi,\bSi) = \bPsi(\bPsi+n_i^{-1}\bSi)^{-1}(\bar \y_i - \bar \X_i^\top\bbe),
\label{eqn:vB}
\end{equation}
where $\yb_i=n_i^{-1}\sum_{j=1}^{n_i}\y_{ij}$.
Thus, we get the estimator
\begin{align}
\btht_a(\bbe, \bPsi,\bSi)=&\c_a^\top\bbe+E[\v_a \mid \y_a]=\c_a^\top \bbe + \vbt_a(\bbe,\bPsi,\bSi)
\non\\
=&
\c_a^\top \bbe + \bPsi(\bPsi+n_a^{-1}\bSi)^{-1}(\yb_a - \Xb_a^\top\bbe),
\label{eqn:Bayes}
\end{align}
which corresponds to the Bayes estimator of $\bth_a$ in the Bayesian framework.

\medskip
When $\bPsi$ and $\bSi$ are known, the maximum likelihood estimator or generalized least squares estimator of $\bbe$ is
\begin{align}
\bbeh(\bPsi,\bSi) =
(\X^\top\D^{-1}\X)^{-1}\X^\top\D^{-1}\y \non,
\label{eqn:beh}
\end{align}
where $\D={\rm block \,diag}(\D_1,\ldots,\D_m)$ and $\D_i=\J_{n_i}\otimes\bPsi+\I_{n_i}\otimes\bSi$ for $i=1,\ldots,m$.
Substituting $\bbeh(\bPsi,\bSi)$ into $\btht_a(\bbe,\bPsi,\bSi)$ yields the estimator
\begin{equation}
\bthh_a(\bPsi,\bSi) = \c_a^\top \hbbe(\bPsi,\bSi) + \bPsi(\bPsi+n_a^{-1}\bSi)^{-1}(\yb_a - \Xb_a^\top\hbbe(\bPsi,\bSi)).
\label{eqn:BLUP}
\end{equation}
It can be easily verified that this estimator is the best linear unbiased predictor (BLUP) of $\bth_a$.

\medskip
We provide consistent estimators of the covariance components $\bSi$ and $\bPsi$.
Concerning  estimation of $\bSi$, it is noted that $E[\{\y_{ij}-\yb_i-(\X_{ij}-\Xb_i)^\top\bbe\}\{\y_{ij}-\yb_i-(\X_{ij}-\Xb_i)^\top\bbe\}^\top]=(1-n_i^{-1})\bSi$ for $i=1, \ldots, m$ and $j=1,\ldots,n_i$, which implies that $\sum_{i=1}^m\sum_{j=1}^{n_i} E[\{\y_{ij}-\yb_i-(\X_{ij}-\Xb_i)^\top\bbe\}\{\y_{ij}-\yb_i-(\X_{ij}-\Xb_i)^\top\bbe\}^\top]=(N-m) \bSi$.
Let $\ybt_i=((\y_{i1}-\yb_i)^\top, \ldots, (\y_{in_i}-\yb_i)^\top)^\top$, $\ybt=(\ybt_1^{\top}, \ldots, \ybt_m^{\top})^\top$, $\Xbt_i=(\X_{i1}-\Xb_i,\ldots,\X_{in_i}-\Xb_i)^\top$ and $\Xbt=(\Xbt_1^{\top}, \ldots, \Xbt_m^{\top})^\top$.
Substituting the statistic $\bbet=(\Xbt^{\top}\Xbt)^{-1}\Xbt^{\top}\ybt$ into $\bbe$, we get an unbised estimator of the form
\begin{equation}
\bSih = {1\over N-m-s_0} (\ybt-\Xbt\bbet)(\ybt-\Xbt\bbet)^\top,
\label{eqn:Sih}
\end{equation}
where $s_0$ is the rank of $\Xbt$.
For estimation of $\bPsi$, it is noted that $E[(\y_{ij}-\X_{ij}^\top\bbe)(\y_{ij}-\X_{ij}^\top\bbe)^\top]=\bPsi+\bSi$ for $i=1, \ldots, m$ and  $j=1,\ldots,n_i$, which implies that $\sum_{i=1}^m \sum_{j=1}^{n_i} E[(\y_i-\X_i\bbe)(\y_i-\X_i\bbe)^\top]=N (\bPsi+\bSi)$.
Substituting the ordinary least squares estimator $\bbeh^{OLS}=(\X^\top\X)^{-1}\X^\top\y$ and $\bSih$ into $\bbe$ and $\bSi$, we get the consistent estimator
\begin{equation}
\bPsih_0 = {1\over N} \sum_{i=1}^m\sum_{j=1}^{n_i} (\y_{ij}-\X_{ij}^\top\bbeh^{OLS})(\y_{ij}-\X_{ij}^\top\bbeh^{OLS})^\top-\bSih.
\label{eqn:Psi0}
\end{equation}
Taking the expectation of $\bPsih_0$, we can see that $E[\bPsih_0]=\bPsi+{\rm Bias}_{\bPsih_0}(\bPsi)$, where
\begin{align}
{\rm Bias}_{\bPsih_0}(\bPsi,\bSi)=&
 {1\over N}\sum_{i=1}^m\sum_{j=1}^{n_i} \X_{ij}^\top(\X^\top\X)^{-1}\X\D^{-1}\X(\X^\top\X)^{-1} \X_{ij}
\non\\
&
-{1\over N}\sum_{i=1}^m\sum_{j=1}^{n_i} (\bSi \X_{ij}^\top-n_i\bPsi\Xb_i^\top) (\X^\top \X)^{-1} \X_{ij}
\non\\
&-{1\over N}\sum_{i=1}^m\sum_{j=1}^{n_i}  \X_{ij}^\top (\X^\top \X)^{-1} (\X_{ij}\bSi-n_i\Xb_i\bPsi),
\label{eqn:PBias}
\end{align}
where $\D={\rm block\ diag}(\D_1, \ldots, \D_m)$ for $\D_i=\J_{n_i}\otimes\bPsi+\I_{n_i}\otimes\bSi$, $i=1, \ldots, m$.
Let $\bPsih_1=\bPsih_0 - {\rm Bias}_{\bPsih_0}(\bPsih_0, \bSih)$.
Then, $\bPsih_1$ is a second-order unbiased estimator of $\bPsi$.
Because $\bPsih_1$ takes a negative value, we modify it as 
\begin{equation}
\bPsih = \H \diag\{ \max(\la_1, 0), \ldots, \max(\la_k, 0)\}\H^\top,
\label{eqn:Psi1}
\end{equation}
where $\H$ is an orthogonal matrix such that $\bPsih_1=\H\diag(\la_1, \ldots, \la_k)\H^\top$.

\medskip
The consistency of $\bSih$ and $\bPsih$ can be shown under the assumptions:

(A1)\ The number of areas $m$ tends to infinity, and $k$, $s$ and $n_i$'s are bounded with respect to $m$.

(A2)\ $\X^\top \X$ is nonsingular and $\X^\top \X/m$ converges to a  positive definite matrix.

\begin{thm}
\label{thm:consistency}
Assume conditions {\rm (A1)} and {\rm (A2)}.
Then, the following asymptotic properties hold for $\bSih$ and $\bPsih$:

{\rm (1)}\ $\bPsih$ is a second-order unbiased estimator of $\bPsi$, while $\bSih$ is an unbiased estimator of $\bSi$.

{\rm (2)}\ $\bSih-\bSi=O_p(m^{-1/2})$, $\bPsih-\bPsi=O_p(m^{-1/2})$ and $\bbeh(\bPsih,\bSih)-\bbe=O_p(m^{-1/2})$.

{\rm (3)}\ For any $\de>0$, $P(\bPsih \not= \bPsih_1)=O(m^{-\de})$.
\end{thm}

The proof is given in the Appendix
Since $\bSih$ and $\bPsih$ are consistent, we can substitute them into (\ref{eqn:BLUP}) to get the empirical best linear unbiased predictor (EBLUP)
\begin{equation}
\bthh_a^{EB} = \c_a^\top \hbbe(\bPsih,\bSih) + \bPsih(\bPsih+n_a^{-1}\bSih)^{-1}(\yb_a - \Xb_a^\top\hbbe(\bPsih,\bSih)).
\label{eqn:EBLUP}
\end{equation}

\section{Evaluation of Uncertainty of EBLUP}
\label{sec:MSE}

The EBLUP suggested in (\ref{eqn:EBLUP}) is expected to have a small estimation error, and it is important to measure how much the estimation error is.
In this section, we derive a second-order approximation of the mean squared error matrix (MSEM) of the EBLUP and provide a second-order unbiased estimator of the MSEM.
The MSEM of the EBLUP is ${\rm MSEM}(\bthh_a^{EB})=E[\{\bthh_a^{EB}-\bth_a\}\{\bthh_a^{EB}-\bth_a\}^\top]$.
It is noted that
$$
\bthh_a^{EB}-\bth_a
=
\{\btht_a(\bbe,\bPsi,\bSi)-\bth_a\} + \{\bthh_a(\bPsi,\bSi)-\btht_a(\bbe,\bPsi,\bSi)\}
+\{\bthh_a^{EB}-\bthh_a(\bPsi,\bSi)\},
$$
where $\btht_a(\bbe,\bPsi,\bSi)$ and $\bthh_a(\bPsi,\bSi)$ are given in (\ref{eqn:Bayes}) and (\ref{eqn:BLUP}).
The following lemma which will proved in the Appendix is useful for evaluating the mean square error matrix.

\begin{lem}
\label{lem:1}
$\bbeh(\bPsi,\bSi)$ is independent of $\y-\X\bbeh^{OLS}$ and $\ybt-\Xbt\bbet$, which implies that $\bbeh(\bPsi,\bSi)$ is independent of $\bSih$ and $\bPsih$.
Also, $\bbeh(\bPsi,\bSi)$ is independent of $\bthh_a^{EB}-\bthh_a(\bPsi,\bSi)$.
\end{lem}

Noting that $\bthh_a(\bPsi,\bSi)-\btht_a(\bbe,\bPsi,\bSi)=\{\c_a^\top-\bPsi(\bPsi+n_a^{-1}\bSi)^{-1}\Xb_a^\top\}\{\bbeh(\bPsi,\bSi)-\bbe\}$, from Lemma \ref{lem:1}, we can decompose the MSEM as
\begin{align}
{\rm MSEM}(\bthh_a^{EB})=&
E[\{\btht_a(\bbe,\bPsi,\bSi)-\bth_a\}\{\btht_a(\bbe,\bPsi,\bSi)-\bth_a\}^\top]
\non\\
&+E[\{\bthh_a(\bPsi,\bSi)-\btht_a(\bbe,\bPsi,\bSi)\}\{\bthh_a(\bPsi,\bSi)-\btht_a(\bbe,\bPsi,\bSi)\}^\top]
\non\\
&+E[\{\bthh_a^{EB}-\bthh_a(\bPsi,\bSi)\}\{\bthh_a^{EB}-\bthh_a(\bPsi,\bSi)\}^\top]
\non\\
=& \G_{1a}(\bPsi,\bSi) + \G_{2a}(\bPsi,\bSi) + E[\{\bthh_a^{EB}-\bthh_a(\bPsi,\bSi)\}\{\bthh_a^{EB}-\bthh_a(\bPsi,\bSi)\}^\top],
\label{eqn:MSE}
\end{align}
where
\begin{equation}
\begin{split}
\G_{1a}(\bPsi,\bSi)=&
(\bPsi^{-1}+n_a\bSi^{-1})^{-1}=n_a^{-1}\bPsi \bLa_a^{-1}\bSi,\\
\G_{2a}(\bPsi,\bSi)=&
(\c_a^\top-\bPsi \bLa_a^{-1}\Xb_a^\top)(\X^\top\D^{-1}\X)^{-1}(\c_a-\Xb_a \bLa_a^{-1}\bPsi),
\end{split}
\label{eqn:G12}
\end{equation}
for $\bLa_a=\bPsi+n_a^{-1}\bSi$.
In the following theorem which will be proved in the Appendix, we approximate the third term as
\begin{align}
\G_{3a}(\bPsi,\bSi)
=&
{n_a^{-2} \over N^2}\bSi\bLa_a^{-1}\sum_{i=1}^m n_i^2\Big\{ \bLa_i\bLa_a^{-1}\bLa_i
+\tr(\bLa_a^{-1}\bLa_i)\bLa_i\Big\}\bLa_a^{-1}\bSi 
\non\\
&+{n_a^{-2} \over N^2(N-m)}(N\bPsi+m\bSi)\bLa_a^{-1}\Big\{\bSi\bLa_a^{-1}\bSi+\tr(\bLa_a^{-1}\bSi)\bSi\Big\}\bLa_a^{-1}(N\bPsi+m\bSi).
\label{eqn:G3}
\end{align}

\begin{thm}
\label{thm:MSE}
The mean squared error matrix of the empirical Bayes estimator $\bthh_a^{EB}$ is approximated as
\begin{equation}
{\rm MSEM}(\bthh_a^{EB})=
\G_{1a}(\bPsi,\bSi) + \G_{2a}(\bPsi,\bSi) + \G_{3a}(\bPsi,\bSi) + O(m^{-3/2}).
\label{eqn:MSE}
\end{equation}
\end{thm}

We next provide a second-order unbiased estimator of the mean squared error matrix of the EBLUP.
A naive estimator of ${\rm MSEM}(\bthh_a^{EB})$ is the plug-in estimator of (\ref{eqn:MSE}) given by $\G_{1a}(\bPsih,\bSih) + \G_{2a}(\bPsih,\bSih) + \G_{3a}(\bPsih,\bSih) $, but this has a second-order bias, because $E[\G_{1a}(\bPsih,\bSih)]=\G_{1a}(\bPsi,\bSi)+O(m^{-1})$.
Correcting this second-order bias, we can derive the second-order unbiased estimator
\begin{equation}
{\rm msem}(\bthh_a^{EB})=\G_{1a}(\bPsih,\bSih) + \G_{2a}(\bPsih,\bSih) + 2\G_{3a}(\bPsih,\bSih).
\label{eqn:msemest}
\end{equation}

\begin{thm}
\label{thm:msemest}
Under the coditions {\rm (A1)} and {\rm (A2)}, it holds that $E[\G_{1a}(\bPsih,\bSih)+\G_{3a}(\bPsih,\bSih)]=\G_{1a}(\bPsi,\bSi)+O(m^{-3/2})$ and 
$$
E[{\rm msem}(\bthh_a^{EB})]={\rm MSEM}(\bthh_a^{EB})+O(m^{-3/2}),
$$
namely, ${\rm msem}(\bthh_a^{EB})$ is a second-order unbiased estimator of ${\rm MSEM}(\bthh_a^{EB})$.
\end{thm}

\section{Confidence Interval for Linear Combination of EBLUP with Corrected Coverage Probability}
\label{sec:ci}

In this section, we consider the confidence interval of the liner combination $\bel^\top\bth_a$ for $\bel \in \mathbb R^k$ for the $a$-th area in the MNER.

\medskip
We begin by estimating the linear combination $\bel^\top\bth_a=\bel^\top(\c_a^\top \bbe+\v_a)$, which is the conditional mean $E[\bel^\top\y_a\mid \v_a]$ given $\v_a$.
A reasonable estimator is provided by the conditional expectation $E[\bel^\top\bth_a\mid \y_a]=\bel^\top\btht_a(\bbe, \bPsi,\bSi)$, where $\btht_a(\bbe, \bPsi,\bSi)$ is given by (\ref{eqn:Bayes}).
By replacing $\bbe$ with the generalized least estimator $\bbeh(\bPsi,\bSi) =(\X^\top\D^{-1}\X)^{-1}\X^\top\D^{-1}\y$, the BLUP of $\bel\bth_a$ is provided by $\bel^\top\bthh_a(\bPsi,\bSi)$, where $\bthh_a(\bPsi,\bSi)$ is given in (\ref{eqn:BLUP}).
Substituting (\ref{eqn:Sih}) and (\ref{eqn:Psi1}) into the BLUP yields the EBLUP $\bel^\top\bthh_a^{EB}$ for $\bthh_a^{EB}$ given in (\ref{eqn:EBLUP}).
The mean squared error is $E[(\bel^\top\bthh_a^{EB}-\bel^\top\bth_a)^2]
= \bel^\top {\rm MSEM}(\bthh_a^{EB})\bel$ and and its second-order unbiased estimator is $\bel^\top{\rm msem}(\bthh_a^{EB})\bel$, namely
$$
E[ \bel^\top{\rm msem}(\bthh_a^{EB})\bel]= E[(\bel^\top\bthh_a^{EB}-\bel^\top\bth_a)^2]+o(m^{-1}),
$$
where ${\rm MSEM}(\bthh_a^{EB})$ and ${\rm msem}(\bthh_a^{EB})$ are given in (\ref{eqn:MSE}) and (\ref{eqn:msemest}).

\medskip
We now construct the confidence interval.
The naive confidence interval is given by 
\begin{equation}
I^{NCI}\ :\ \bel^\top \bthh_a^{EB}\pm z_{\al/2}\times\sqrt{\mathstrut \bel^\top{\rm msem}(\bthh_a^{EB}) \bel},
\label{eqn:CIn}
\end{equation}
where $z_{\al/2}$ is the $100(1-\al/2)\%$ percentile of the standard normal distribution.
However, this confidence interval does not have the second-order accuracy, namely $P(\bel^\top\bth_a \in I^{NCI})=1-\al+O(m^{-1})$.
To derive a confidence interval with second-order accuracy, we need to evaluate the second moment of $\bel^\top{\rm msem}(\bthh_a^{EB})\bel-\bel^\top {\rm MSEM}(\bthh_a^{EB})\bel$.
Let
\begin{align}
\begin{split}
V(\bthh_a^{EB})=&{n_a^{-4} \over N^2}\sum_{i=1}^m n_i^2 \Big\{ (\bel^\top\bSi\bLa_a^{-1}\bLa_i\bel)^2+\bel^\top\bSi\bLa_a^{-1}\bLa_i\bLa_a^{-1} \bSi\bel \times \bel^\top\bLa_i\bel \Big\}
\\
&+{2n_a^{-4} m^2\over N^2(N-m)}(\bel^\top\bSi\bLa_a^{-1}\bSi\bLa_a^{-1} \bSi\bel)^2+{2n_a^{-2} \over N-m}(\bel^\top\bPsi\bLa_a^{-1}\bSi\bLa_a^{-1} \bPsi\bel)^2
\\
&-{2n_a^{-3} m\over N(N-m)}\Big\{ \bel^\top\bPsi\bLa_a^{-1}\bSi\bLa_a^{-1} \bPsi\bel \times \bel^\top\bSi\bLa_a^{-1}\bSi\bel+(\bel^\top\bSi\bLa_a^{-1}\bSi\bLa_a^{-1}\bPsi\bel)^2 \Big\}.
\end{split}
\label{eqn:vmse}
\end{align}

\begin{lem}
\label{lem:ci}
Under the coditions {\rm (A1)} and {\rm (A2)}, it holds that
\begin{align*}
E\Big[\Big\{\bel^\top{\rm msem}(\bthh_a^{EB})\bel-\bel^\top {\rm MSEM}(\bthh_a^{EB})\bel \Big\}^2\Big]=V(\bthh_a^{EB})+o(m^{-1}),
\end{align*}
and for $c\geq 3$,
\begin{align*}
E\Big[\Big\{\bel^\top{\rm msem}(\bthh_a^{EB})\bel-\bel^\top {\rm MSEM}(\bthh_a^{EB})\bel \Big\}^c\Big]=o(m^{-1}).
\end{align*}
\end{lem}

\begin{thm}
\label{thm:ci}
Under the coditions {\rm (A1)} and {\rm (A2)}, it holds that for any $z$, 
\begin{align*}
P\Big(& \bel^\top\bthh_a^{EB}-\{\bel^\top{\rm msem}(\bthh_a^{EB})\bel\}^{1/2}z \leq \bel^\top\bth_a \leq \bel^\top\bthh_a^{EB}+\{\bel^\top{\rm msem}(\bthh_a^{EB})\bel\}^{1/2}z \Big)
\\
=&
2\Phi(z)-1- { V(\bthh_a^{EB})\over 4\{\bel^\top {\rm MSEM}(\bthh_a^{EB}) \bel\}^2 } (z^3+z)\phi(z)+o(m^{-1}),
\end{align*}
where $\Phi(\cdot)$ and $\phi(\cdot)$ are the distribution and density functions of the standard normal distribution.
\end{thm}

Solving the equation
$$
2\Phi(z)-1- { V(\bthh_a^{EB})\over 4\{\bel^\top {\rm MSEM}(\bthh_a^{EB}) \bel\}^2 } (z^3+z)\phi(z)=1-\al,
$$
we get the solution given by 
$$
z^*=z_{\alpha/2}+(z_{\alpha/2}^3+z_{\alpha/2}) V(\bthh_a^{EB}) / 8\{\bel^\top {\rm MSEM}(\bthh_a^{EB}) \bel\}^2,
$$
which provides the improved confidence interval
\begin{equation}
I^{ICI}\ :\ \bel^\top\bthh_a^{EB} \pm \{\bel^\top{\rm msem}(\bthh_a^{EB})\bel\}^{1/2}z^*.
\label{eqn:ICI}
\end{equation}
Then from Theorem \ref{thm:ci}, it follows that $P(\bel^\top\bth_a\in I^{ICI})=1-\al+o(m^{-1})$.

\section{Simulation and Empirical  Studies}
\label{sec:sim}

\subsection{Finite sample performances}

We now investigate finite sample performances of EBLUP in terms of MSEM and the second-order unbiased estimator of MSEM by simulation.

\bigskip
{\bf [1] \ Setup of simulation experiments}.\ \ 
We treat the multivariate Nested-Error model,
\begin{equation}
\y_{ij} = \X_{ij}^\top \bbe + \v_i + \bep_{ij},
\quad i=1, \ldots, m, 
\quad j=1, \ldots, n_i.
\non
\end{equation}
We take $m=40$, $k=2,3$ and $\bbe=(0.8,-0.5,-0.3,0.6)^\top$ for $k=2$ and $\bbe=(0.8,-0.5,-0.3,0.6,0.4,-0.2)^\top$ for $k=3$.
Moreover, we equally divided areas into four groups $(G = 1,\ldots,4)$, so that each group has ten areas and the areas in the same group has the same sample size $n_G = 3G-2$. 
The design matrix, $\X_{ij}$ is $2k \times k$ matrix, such that 
$$
\X_{ij}=\begin{pmatrix} 1 & x_{i1} & 0 & 0\\ 0&0&1 & x_{i2} \end{pmatrix}^\top, \X_{ij}=\begin{pmatrix} 1 & x_{i1} & 0 & 0 & 0 & 0\\ 0&0&1 & x_{i2} & 0 & 0\\ 0&0&0&0& 1& x_{i3} \end{pmatrix}^\top
$$
for $k=2,3$ respectevely.
We generate $x_{ij}$ from uniform distribution on $(-1,1)$, which are fixed through the simulation runs.
As a setup of the covariance matrix $\bPsi$ of the random effects, we consider 
$$
\bPsi= \left\{\begin{array}{ll}
\rho\bpsi_2\bpsi_2^\top+(1-\rho){\rm diag}(\bpsi_2\bpsi_2^\top) & {\rm for}\ k=2,\\
\rho\bpsi_3\bpsi_3^\top+(1-\rho){\rm diag}(\bpsi_3\bpsi_3^\top) & {\rm for}\ k=3,
\end{array}\right.
$$
where $\bpsi_2=(\sqrt{1.5}, \sqrt{0.5})^\top$, $\bpsi_3=(\sqrt{1.5}, 1, \sqrt{0.5})^\top$, and ${\rm diag}(\A)$ denotes the diagonal matrix consisting of diagonal elements of matrix $\A$.
Here, $\rho$ is the correlation coefficient, and  we handle the three cases $\rho=0.25, 0.5, 0.75$.
The cases of negative correlations are omitted, because we observe the same results with those of positive ones.
Concerning the dispersion matrices $\bSi$ of sampling errors $\bep_i$, we set $\bSi=\I_{k}$.
We consider three patterns of distribution of $\v_i$, that is, M1: $\v_i$ is normally distributed, M2: $\v_i$ follows multivariate $t$ distribution with degrees of freedom $5$ and M3: $\v_i$ follows multivariate chi-squared distribution with degrees of freedom $2$.
The distribution of $\ep_i$ is normal.

\begin{figure}[htbp]
  \begin{center}
    \begin{tabular}{ccc}
      \begin{minipage}{0.33\hsize}
        \begin{center}
          \includegraphics[clip, width=5.5cm]{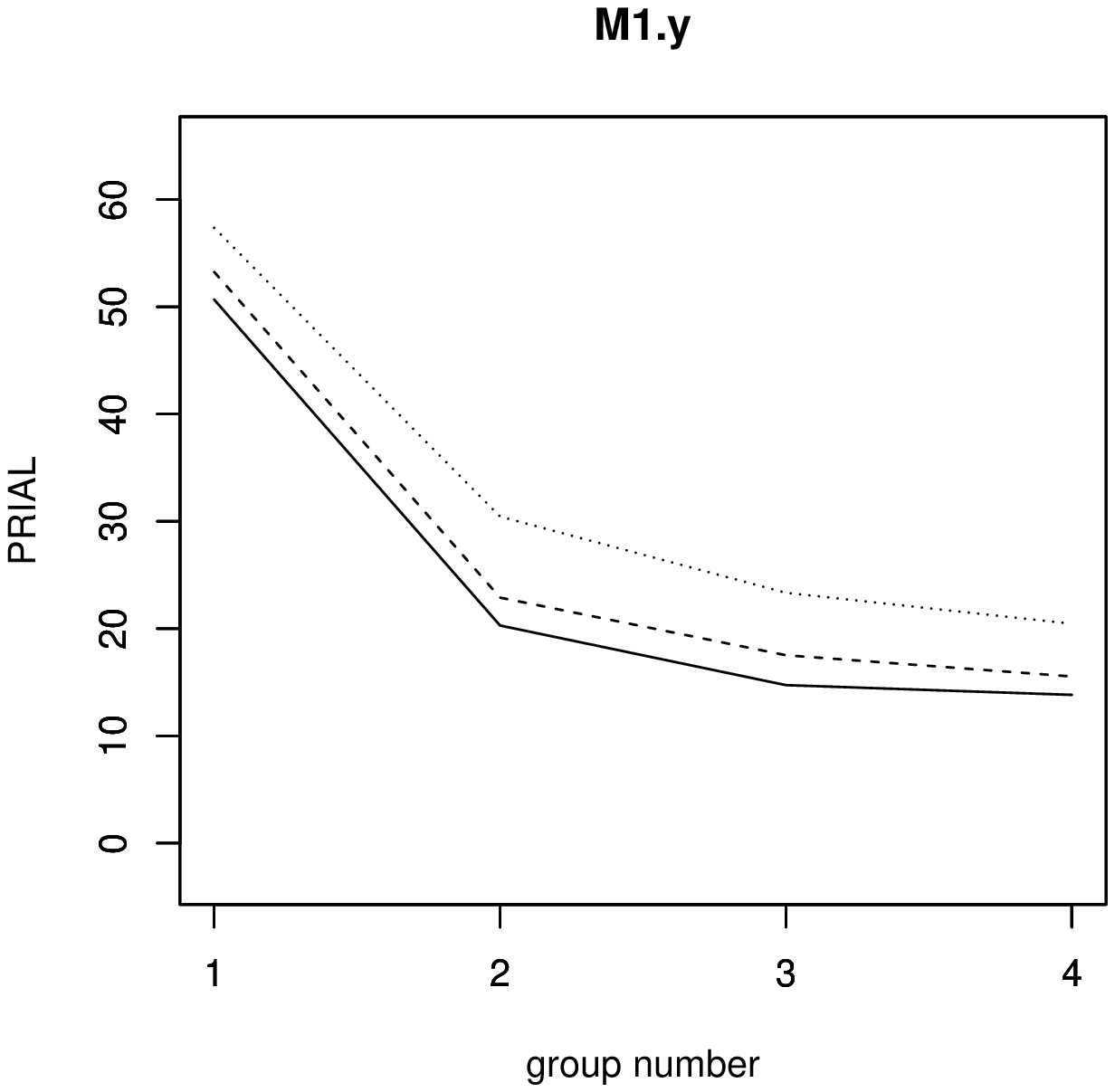}
          \hspace{1.6cm} 
        \end{center}
      \end{minipage}
      \begin{minipage}{0.33\hsize}
        \begin{center}
          \includegraphics[clip, width=5.5cm]{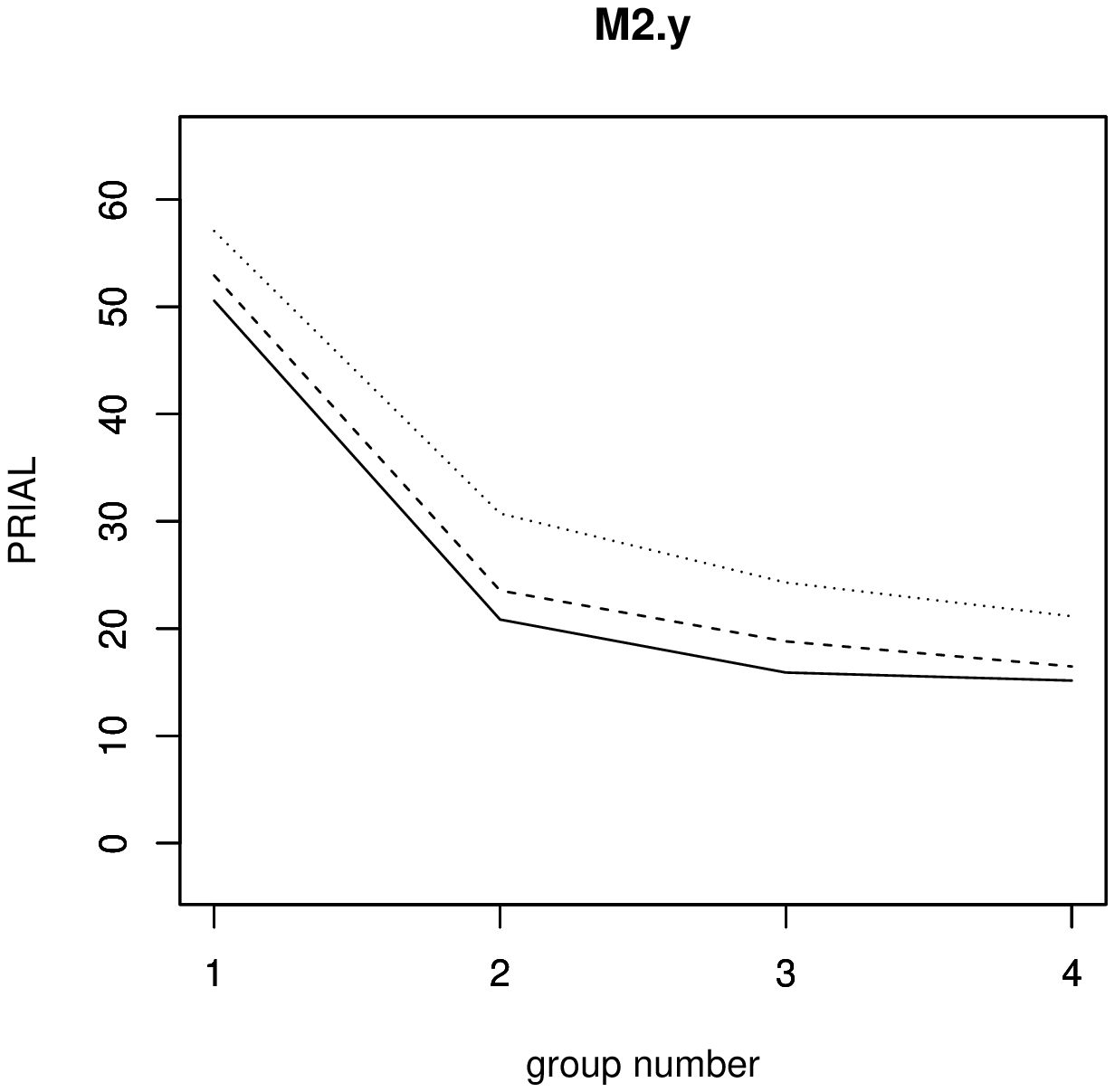}
          \hspace{1.6cm} 
        \end{center}
      \end{minipage}
      \begin{minipage}{0.33\hsize}
        \begin{center}
          \includegraphics[clip, width=5.5cm]{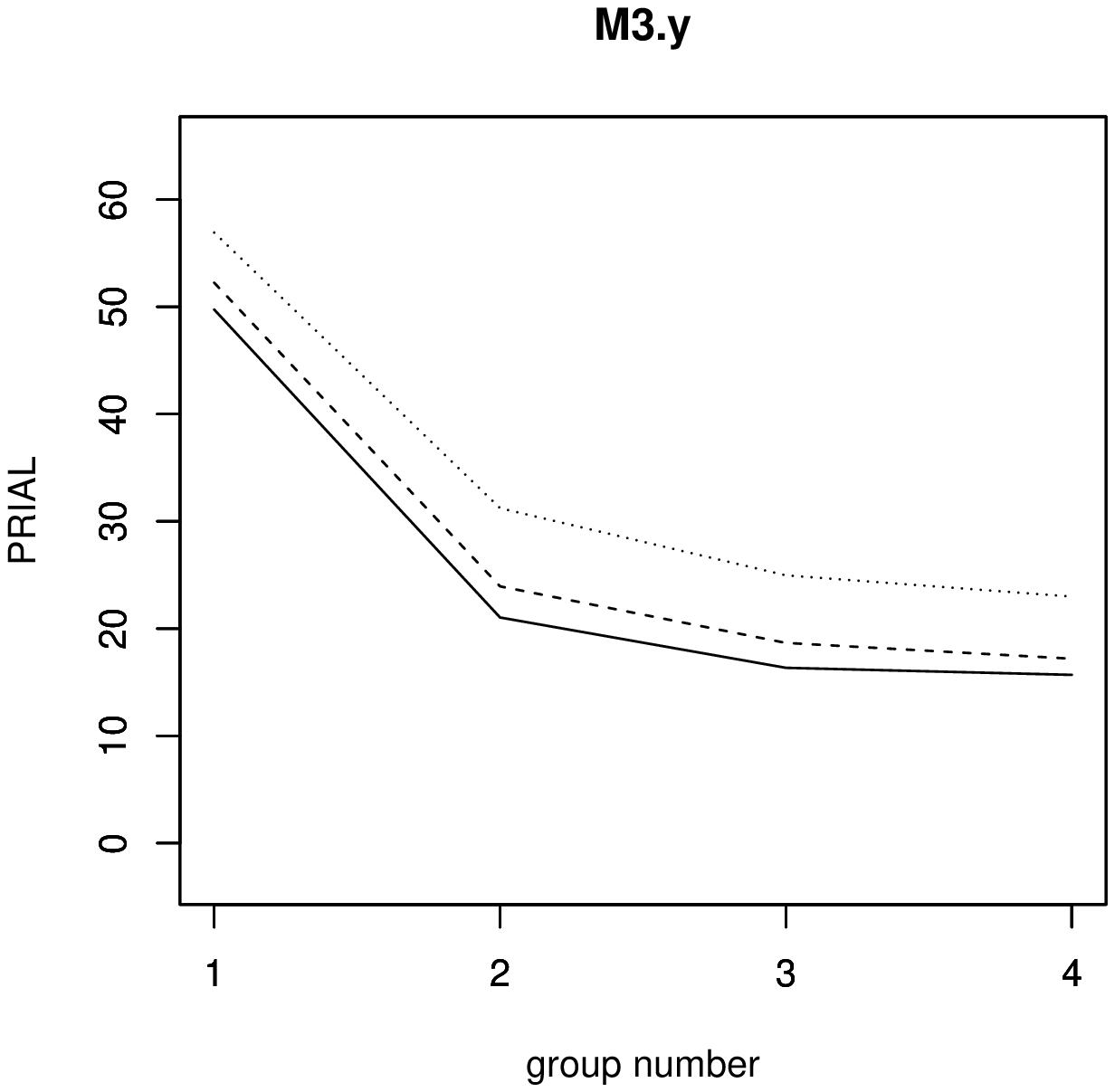}
          \hspace{1.6cm} 
        \end{center}
      \end{minipage}\\
      \begin{minipage}{0.33\hsize}
        \begin{center}
          \includegraphics[clip, width=5.5cm]{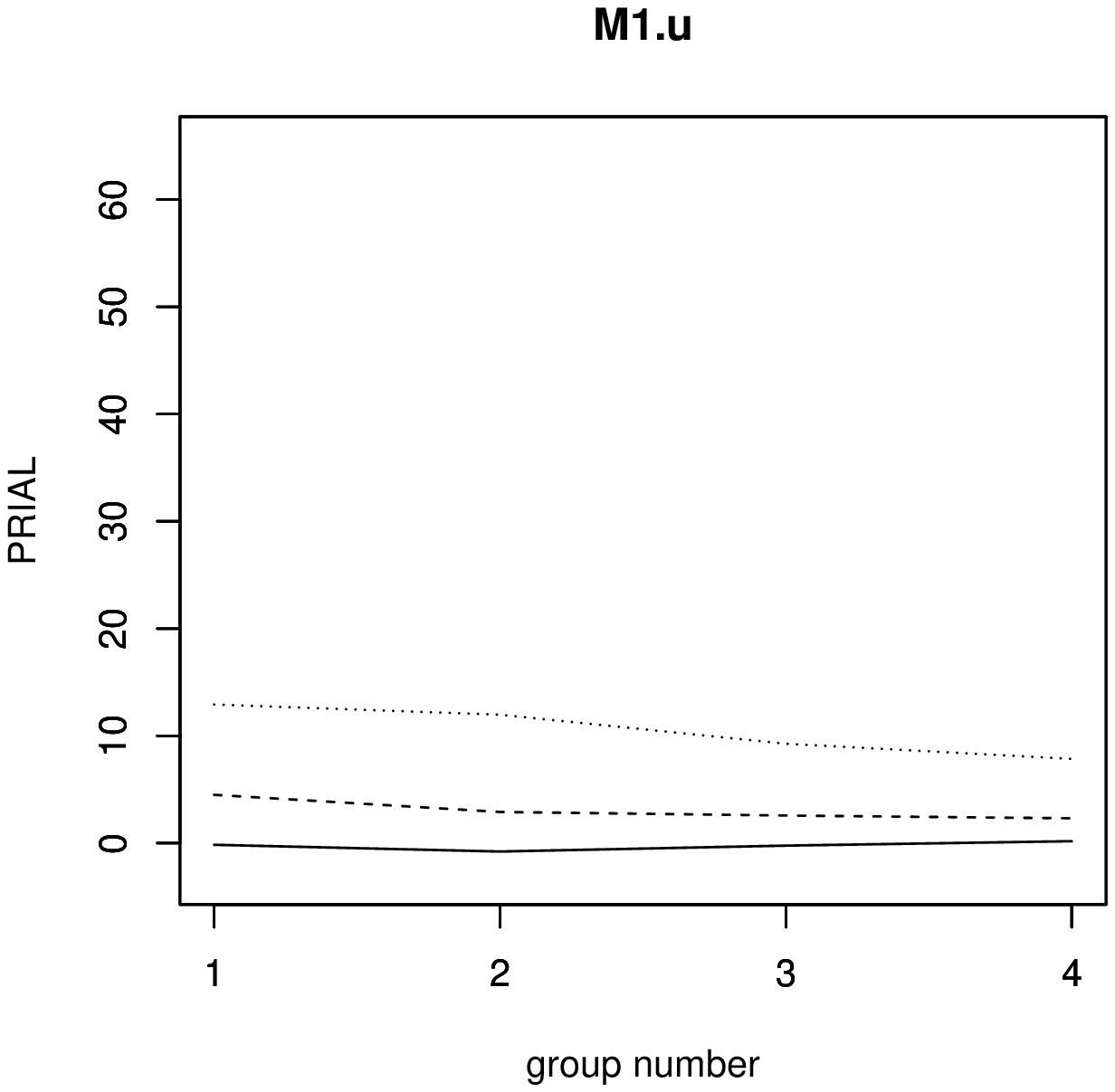}
          \hspace{1.6cm} 
        \end{center}
      \end{minipage}
      \begin{minipage}{0.33\hsize}
        \begin{center}
          \includegraphics[clip, width=5.5cm]{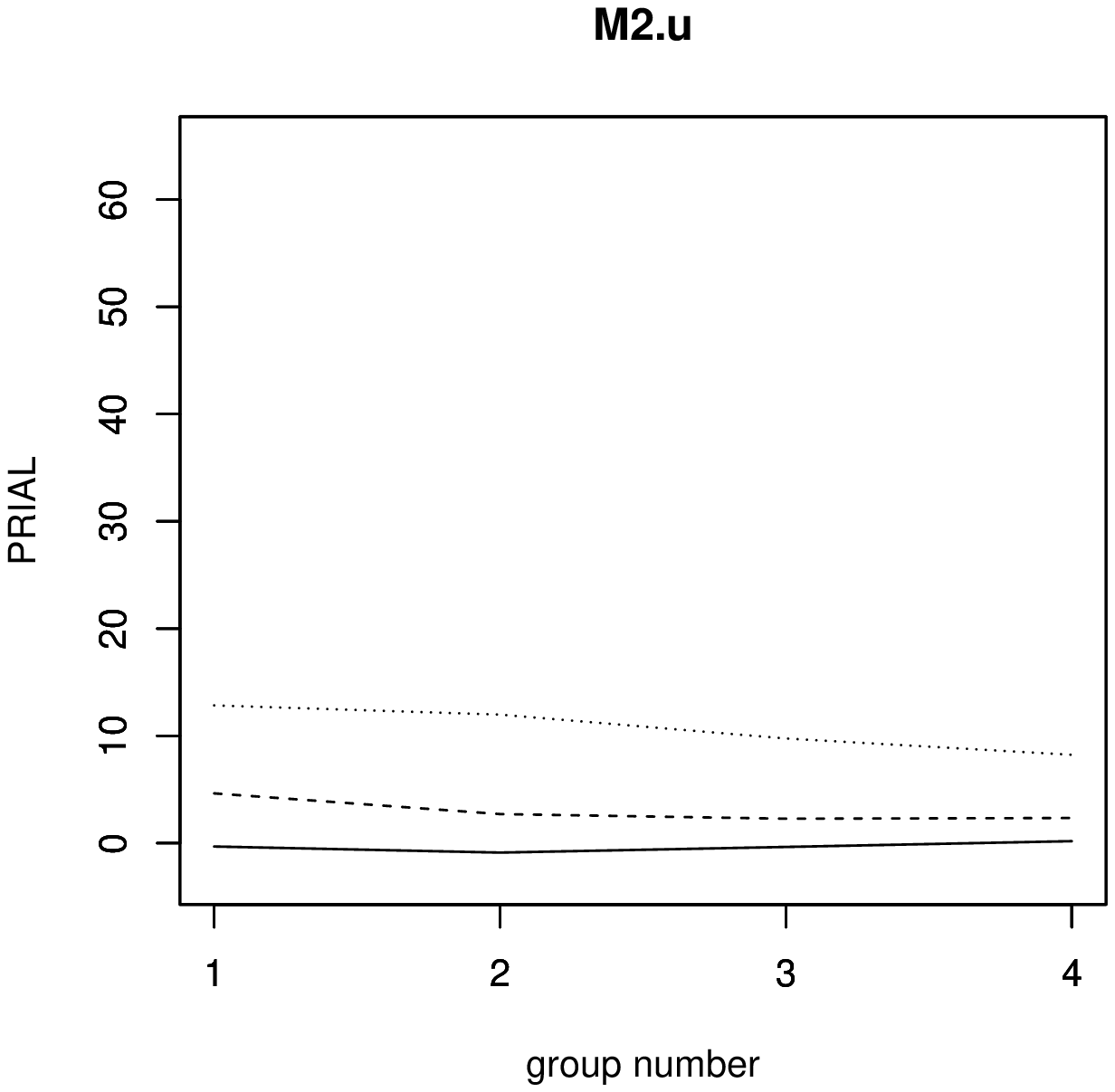}
          \hspace{1.6cm} 
        \end{center}
      \end{minipage}
      \begin{minipage}{0.33\hsize}
        \begin{center}
          \includegraphics[clip, width=5.5cm]{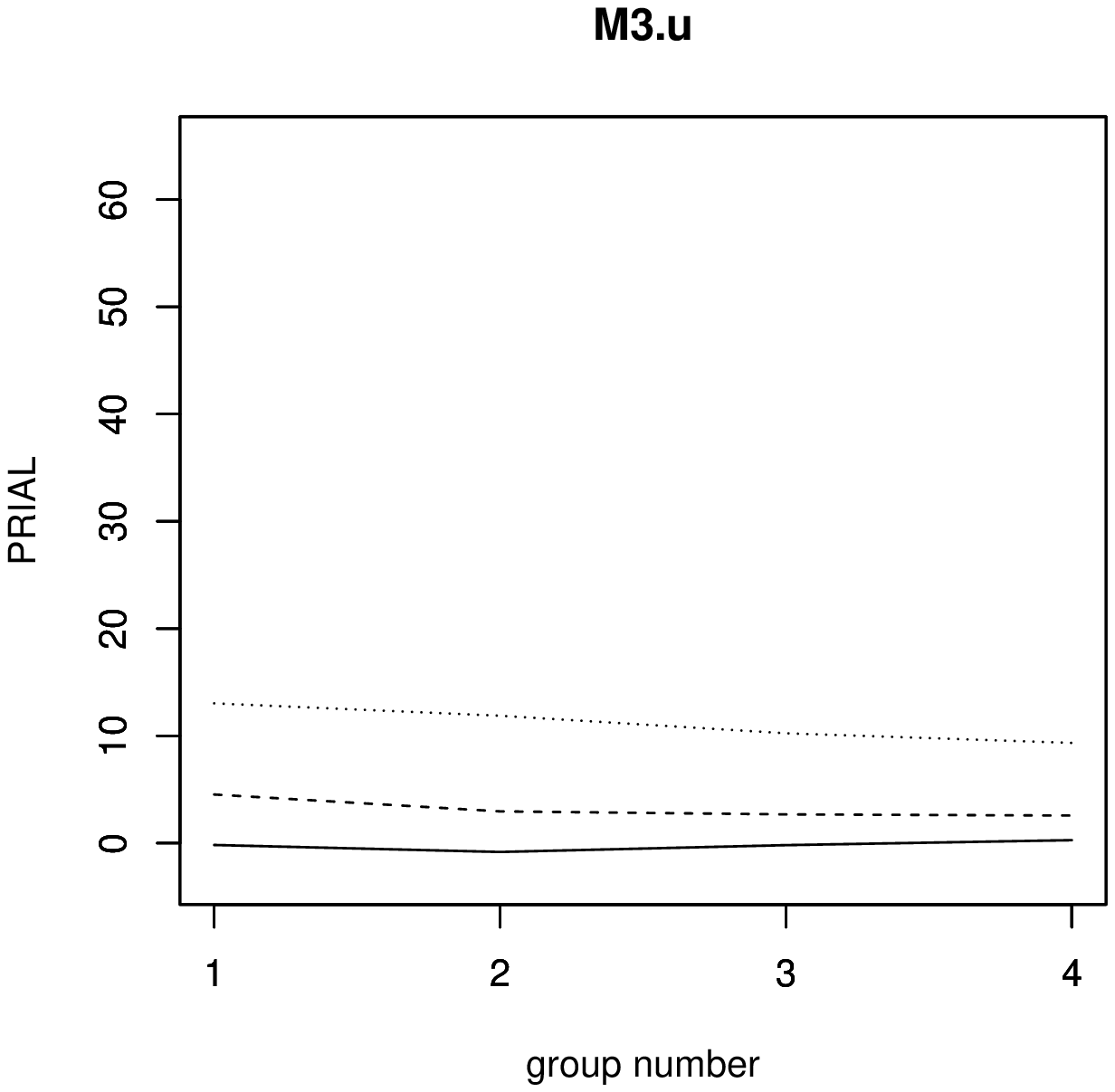}
          \hspace{1.6cm} 
        \end{center}
      \end{minipage}
    \end{tabular}
    \caption{PRIAL for $\rho=0.25$ (real line), $\rho=0.5$ (dashed line) and $\rho=0.75$ (dotted line) in case of $k=2$.}
    \label{fig:prial}
  \end{center}
\end{figure}

\begin{figure}[htbp]
  \begin{center}
    \begin{tabular}{ccc}
      \begin{minipage}{0.33\hsize}
        \begin{center}
          \includegraphics[clip, width=5.5cm]{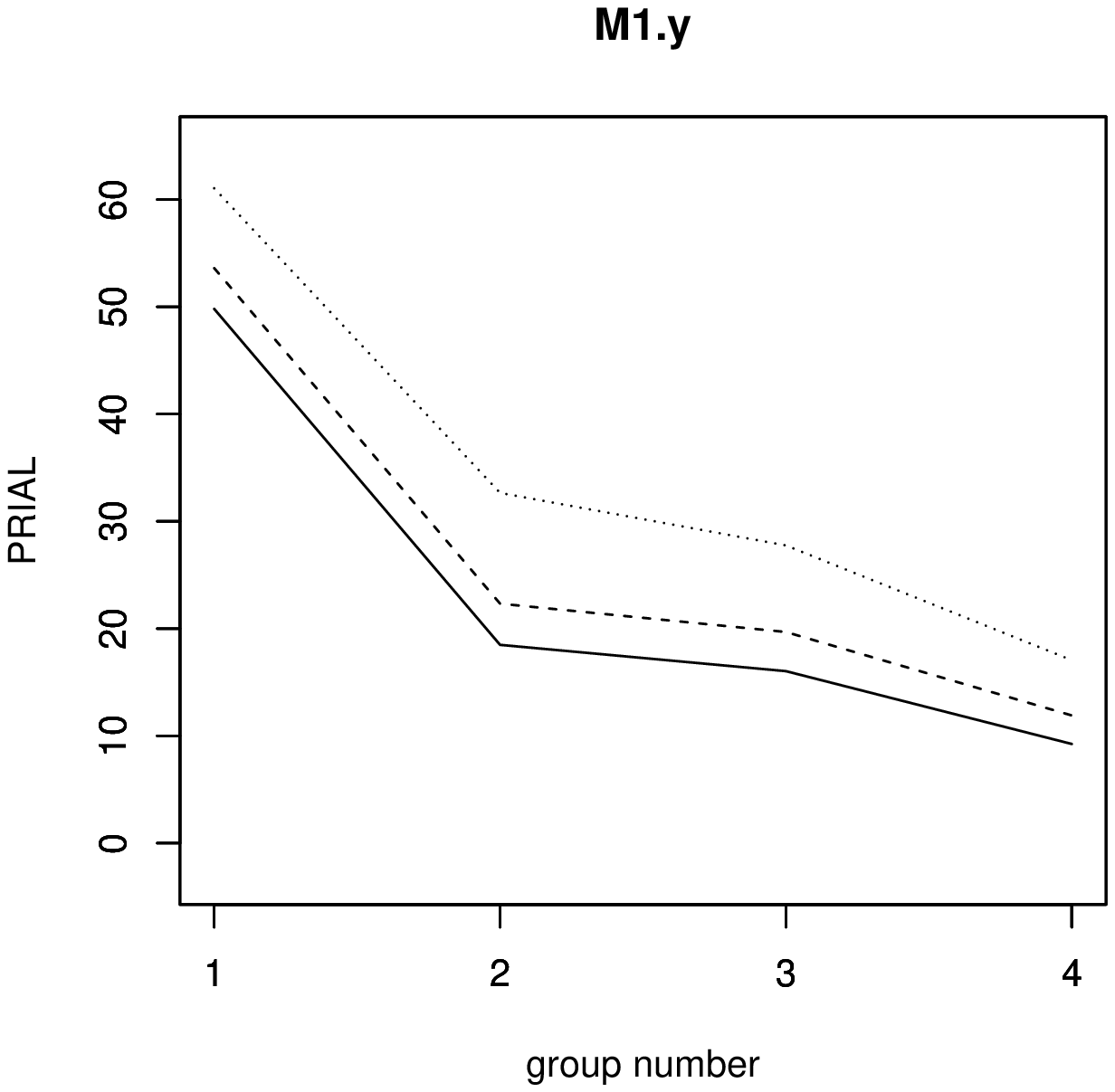}
          \hspace{1.6cm} 
        \end{center}
      \end{minipage}
      \begin{minipage}{0.33\hsize}
        \begin{center}
          \includegraphics[clip, width=5.5cm]{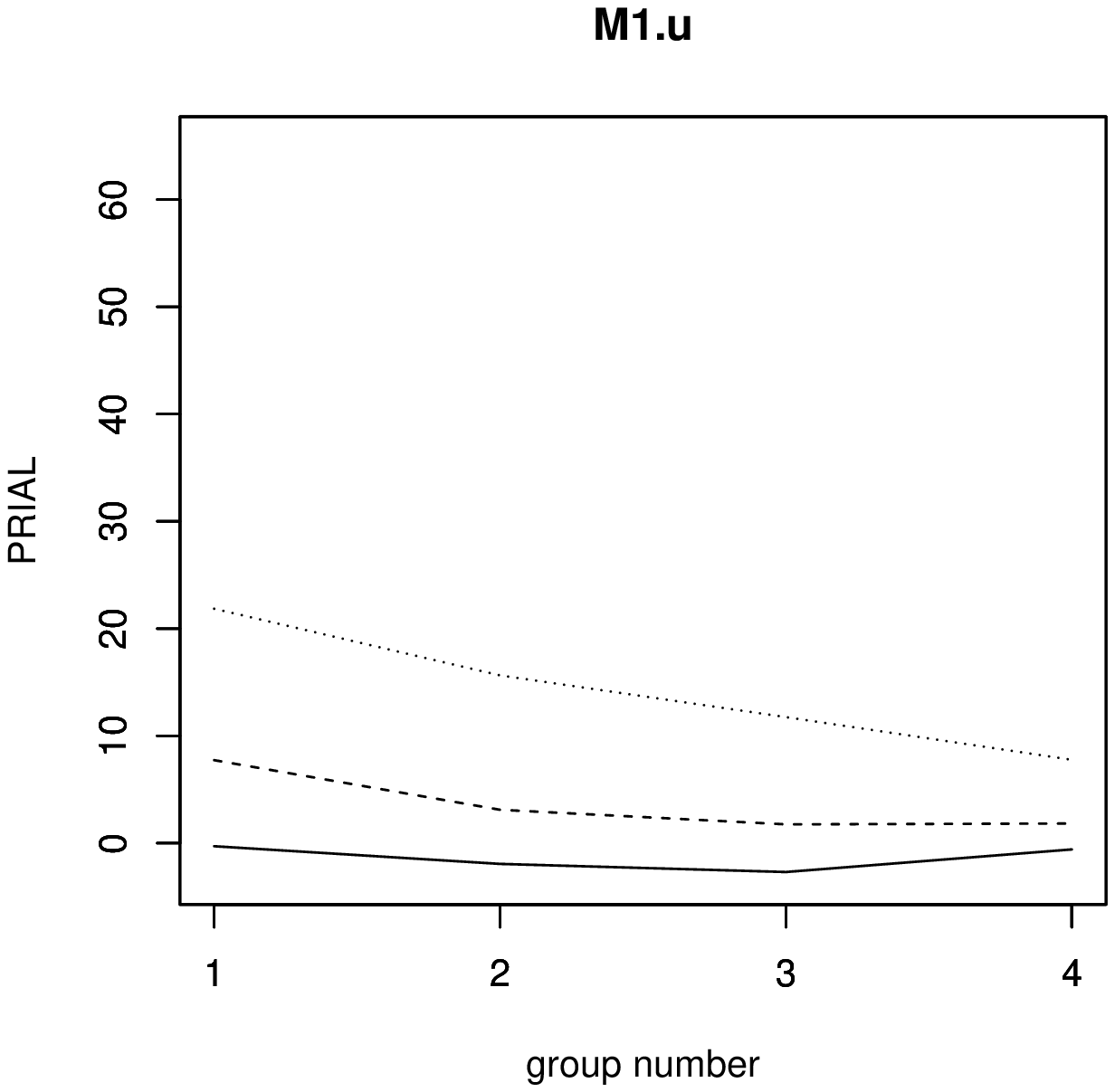}
          \hspace{1.6cm} 
        \end{center}
      \end{minipage}
    \end{tabular}
    \caption{PRIAL for $\rho=0.25$ (real line), $\rho=0.5$ (dashed line) and $\rho=0.75$ (dotted line) in case of $k=3$.}
    \label{fig:prialk3}
  \end{center}
\end{figure}

\bigskip
{\bf [2]\ Comparison of MSEM}.\ \ 
We begin with obtaining the true mean squared error matrices of the EBLUP $\bthh_a^{EB} = \bthh_a(\bPsih,\bSih)$ by simulation.
Let $\{\y_i^{(r)}, i=1, \ldots, m\}$ be the simulated data in the $r$-th replication for $r=1,\ldots, R$ with $R=50,000$.
Let $\bPsih^{(r)}$, $\bSih^{(r)}$ and $\bth_a^{(r)}$ be the values of $\bPsih$, $\bSih$ and $\bth_a=\Xb_{a}^\top \bbe + \v_a$ in the $r$-th replication.
Then the simulated value of the true mean squared error matrices is calculated by
\begin{align}
{\rm MSEM}(\bthh_a^{EB}) = R^{-1}\sum_{i=1}^R \big\{\bthh_a(\bPsih^{(r)},\bSih^{(r)})-\bth_a^{(r)}\big\}\big\{\bthh_a(\bPsih^{(r)},\bPsih^{(r)})-\bth_a^{(r)}\big\}^\top.
\label{eqn:simmsem}
\end{align}
To measure relative improvement of EBLUP, we calculate the percentage relative improvement in the average loss (PRIAL) of $\bthh_a^{EB}$ over $\y_a$, defined by
$$
{\rm PRIAL}(\bthh_a^{EB}, \y_a )
= 100 \times \Big[ 1 - { \tr\{ {\rm MSEM}(\bthh_a^{EB})\} \over \tr\{ {\rm MSEM}(\y_a)\} } \Big].
$$
It is also interesting to compare $\bthh_a^{EB}$ with the EBLUP $\bthh_a^{uEB}$ derived from the univariate Nestd-Error model.
Thus, we calculate the PRIAL given by
$$
{\rm PRIAL}(\bthh_a^{EB}, \bthh_a^{uEB} )
= 100 \times \Big[ 1 - { \tr\{ {\rm MSEM}(\bthh_a^{EB})\} \over \tr\{ {\rm MSEM}(\bthh_a^{uEB})\} } \Big],
$$
and those values are reported in Figure \ref{fig:prial} and \ref{fig:prialk3}.

\medskip
Figure \ref{fig:prial} reports the PRIAL for $k=2$ and three patterns of distribution of $\v_i$; M1, M2 and M3.
We can see that the performances of $\bthh_a^{EB}$ are stable regardless of the distribution of $\v_i$.
In all the cases, $\bthh_a^{EB}$ improves on $\y_a$ largely and the improvement rates are larger for larger $\rho$; for normal case (M1).
On the other hand, $\bthh_a^{EB}$ improves on $\bthh_a^{uEB}$ for large $\rho$, but the univariate EBLUP $\bthh_a^{uEB}$ is slightly better than $\bthh_a^{EB}$ for $\rho=0.25$ for some areas, but the difference is not significant.
This is because the low accuracy in estimation of the covariance matrix $\bPsi$ and $\bSi$ has more adverse influence on prediction than the benefit from incorporating the small correlation into the estimation.
Moreover, the PRIAL is larger for the groups with small sample size.
This is reasonable because the benefit given by incorporating the information from neibouring areas is large for such groups.

Figure \ref{fig:prialk3} reports the PRIAL for $k=3$ and a pattern of distribution of $\v_i$; M1.
The results are almost the same with the case for $k=2$.
The PRIAL is larger for $k = 3$ than for $k = 2$ in the case of $\rho=0.75$, but smaller in the case of $\rho=0.25$.
This is because when $m$ is fixed as $m = 40$, the accuracy in estimation of the covariance matrices gets smaller for
the larger dimension.

\bigskip
{\bf [3]\ Finite sample performances of the MSEM estimator}.\ \ 
We next investigate the performance of the second-order unbiased estimator ${\rm msem}(\bthh_a^{EB})$ of MSEM given in Theorem \ref{thm:msemest}.
We use the same data generating process as mentioned above and we take only $k=2$.
We consider the normal case (M1) as a pattern of distributions for $\v_i$.
The simulated values of the MSEM are obtained from (\ref{eqn:simmsem}) based on $R = 50,000$ simulation runs.
Then, based on $R = 5,000$ simulation runs, we calculate the relative bias (RB) of MSEM estimators given by
\begin{align*}
{\rm RB}_a={1 \over R} \sum_{r=1}^R {{\rm msem}(\bthh_a^{EB(r)})-{\rm MSEM}(\bthh_a^{EB}) \over {\rm MSEM}(\bthh_a^{EB})}
\end{align*}
where ${\rm msem}(\bthh_a^{EB(r)})$ is the MSEM estimator in the $r$-th replication.
In Table \ref{tab:rb}, we report mean values of ${\rm RB}_a$ in each group.
For comparison, results for the naive MSEM estimator, without any bias correction, are reported in Table \ref{tab:rb} as well.
The naive MSEM estimator is the plug-in estimator of the asymptotic MSEM (\ref{eqn:MSE}).
The relative bias is small for the diagonal elements, less than $10\%$ in almost the cases, whereas considerably large for off-diagonal elements.
The naive MSEM estimator is more biased than the analytical MSEM estimator for diagonal elements in all cases, so that the bias correction in MSEM estimator is successful.
On the other hand, the analytical MSEM estimator is more biased slightely than the naive MSEM estimator for off-diagonal elements in some cases.

\begin{table}[htbp]
\begin{center}
\resizebox{16cm}{!} {
\begin{tabular}{ccccccc}
 \hline
&\multicolumn{2}{c}{$\rho=0.25$}&\multicolumn{2}{c}{$\rho=0.5$}&\multicolumn{2}{c}{$\rho=0.75$} \\
 & RB & NRB & RB & NRB &RB & NRB  \\
  \hline
$G_1$ & $
\left[
\begin{array}{rr}
-3.5 & 0.8\\
0.8 & -8.1
\end{array}
\right]
$
&
$
\left[
\begin{array}{rr}
-6.5 & -0.6\\
-0.6 & -12.2
\end{array}
\right]
$
&
$
\left[
\begin{array}{rr}
-4.6 & -1.8\\
-1.8 & -7.7
\end{array}
\right]
$
&
$
\left[
\begin{array}{rr}
-7.7 & -3.2\\
-3.2 & -11.8
\end{array}
\right]
$
&
$
\left[
\begin{array}{rr}
-3.5 & -2.7\\
-2.7 & -7.6
\end{array}
\right]
$
&
$
\left[
\begin{array}{rrr}
-6.3 & -4.3\\
-4.3 & -11.7
\end{array}
\right]
$
\\

$G_2$ & $
\left[
\begin{array}{rr}
-1.9 & 26.7\\
26.7 & -4.1
\end{array}
\right]
$
&
$
\left[
\begin{array}{rr}
-4.0 & 36.1\\
36.1 & -8.8
\end{array}
\right]
$
&
$
\left[
\begin{array}{rr}
-0.5 & 0.9\\
0.9 & -3.1
\end{array}
\right]
$
&
$
\left[
\begin{array}{rr}
-2.8 & 8.1\\
8.1 & -8.5
\end{array}
\right]
$
&
$
\left[
\begin{array}{rr}
-0.6 & 3.7\\
3.7 & -4.7
\end{array}
\right]
$
&
$
\left[
\begin{array}{rrr}
-3.3 & 9.8\\
9.8 & -11.4
\end{array}
\right]
$
\\

$G_3$ & $
\left[
\begin{array}{rr}
-0.7 & 54.8\\
54.8 & -2.7
\end{array}
\right]
$
&
$
\left[
\begin{array}{rr}
-2.3 & 77.1\\
77.1 & -6.9
\end{array}
\right]
$
&
$
\left[
\begin{array}{rr}
-0.4 & 19.3\\
19.3 & -2.1
\end{array}
\right]
$
&
$
\left[
\begin{array}{rr}
-2.3 & 37.5\\
37.5 & -7.2
\end{array}
\right]
$
&
$
\left[
\begin{array}{rr}
-1.7 & 9.9\\
9.9 & -2.9
\end{array}
\right]
$
&
$
\left[
\begin{array}{rrr}
-4.3 & 26.9\\
26.9 & -10.6
\end{array}
\right]
$
\\

$G_4$ & $
\left[
\begin{array}{rr}
-0.2 & 30.2\\
30.2 & -0.5
\end{array}
\right]
$
&
$
\left[
\begin{array}{rr}
-1.5 & 56.5\\
56.5 & -4.2
\end{array}
\right]
$
&
$
\left[
\begin{array}{rr}
0.5 & 5.6\\
5.6 & -0.3
\end{array}
\right]
$
&
$
\left[
\begin{array}{rr}
-1.1 & 29.9\\
29.9 & -5.0
\end{array}
\right]
$
&
$
\left[
\begin{array}{rr}
0.4 & 0.1\\
0.1 & 0.6
\end{array}
\right]
$
&
$
\left[
\begin{array}{rrr}
-2.3 & 28.1\\
28.1 & -7.8
\end{array}
\right]
$
\\
\hline
\end{tabular}
}
\caption{{The Mean Values of Percentage Relative Bias in Each Group (RB) and Relative Bias of Naive MSE Estimator (RBN).}}
\label{tab:rb}
\end{center}
\end{table}

\bigskip
{\bf [4]\ Finite sample performances of the confidence interval}.\ \ 
We investigate the performance of the improved confidence interval given in (\ref{eqn:ICI}).
Table \ref{tab:ci} reports values of coverage probabilities (CP) and average length (AL) for $1-\al=95\%$ confidence coefficient, where the setup of the simulation experiment is the same as above, namely the three patterns of distributions of $\v_i$, M1, M2 and M3 and the three cases of $\rho=0.25, 0.5, 0.75$ are treated.
Table \ref{tab:ci} also reports values of  CP and AL in parentheses for the naive confidence interval (\ref{eqn:CIn}).

\medskip
For all patterns of distributions of $\v_i$ and correlation coefficients, values of CP are close to the nominal level of $0.95$ and are higher than those for the naive method, especially for areas with small sample size.
This is coincident with Diao et al. (2014), which considered the confidence interval estimator under the Fay-Herriot model. 
Values of CP for areas with large sample sizes are slightly higher than those for the naive method, but the differences are negligibly small.
Values of AL are also larger than those for the naive method for areas with small sample size, and the diferrence is negligible for areas with large sample size.

\begin{table}[htbp]
\begin{center}
\resizebox{16cm}{!} {
\begin{tabular}{ccccccccccc}
 \hline
&&\multicolumn{3}{c}{Normal}&\multicolumn{3}{c}{t}&\multicolumn{3}{c}{chi-square} \\
\cline{3-5}\cline{6-8}\cline{9-11}
$\rho$& & 0.25 & 0.5 & 0.75 & 0.25 & 0.5 & 0.75 & 0.25 & 0.5 & 0.75  \\
\hline
$G_1$& CP & 0.949 & 0.950 & 0.957 & 0.946 & 0.945 & 0.950 & 0.940 & 0.940 & 0.943\\ 
&            & (0.936) & (0.933) & (0.928) & (0.935) & (0.929) & (0.923) & (0.927) & (0.924) & (0.918)\\ 
& AL & 3.700 & 3.330 & 2.768 & 3.655 & 3.279 & 2.733 & 3.615 & 3.239 & 2.705   \\
&      & (3.501) & (3.108) & (2.512) & (3.459) & (3.073) & (2.472) & (3.412) & (3.023) & (2.434)   \\
$G_2$& CP & 0.947 & 0.943 & 0.941 & 0.945 & 0.941 & 0.937 & 0.941 & 0.940 &0.935 \\  
&               & (0.945) & (0.940) & (0.937) & (0.943) & (0.940) & (0.932) & (0.939) & (0.938) & (0.931) \\
& AL & 2.340 & 2.239 & 1.978 & 2.357 & 2.218 & 1.948 & 2.339 & 2.197 & 1.924  \\
&      & (2.364) & (2.226) & (1.950) & (2.343) & (2.204) & (1.925) & (2.321) & (2.177) & (1.898)  \\
$G_3$& CP & 0.947 & 0.947 & 0.946 & 0.947 & 0.947 & 0.947 & 0.948 & 0.946 & 0.945\\ 
&               & (0.946) & (0.946) & (0.944) & (0.947) & (0.946) & (0.945) & (0.947) & (0.944) & (0.943)\\  
& AL & 1.920 & 1.845 & 1.706 & 1.905 & 1.833 & 1.687 & 1.896 &1.822 & 1.673  \\
&      & (1.909) & (1.839) & (1.690) & (1.897) & (1.826) & (1.676) & (1.886) & (1.811) & (1.661)  \\
$G_4$& CP & 0.949 & 0.949 & 0.951 & 0.951 & 0.949 & 0.951 & 0.950 & 0.951 & 0.952\\  
&               & (0.947) & (0.948) & (0.949) & (0.950) & (0.948) & (0.950) & (0.949) & (0.950) & (0.950)\\ 
& AL & 1.653 &1.607 & 1.531 & 1.643 & 1.600 & 1.520 & 1.639 & 1.595 & 1.514   \\
&      & (1.644) &(1.602) & (1.520) & (1.636) & (1.594) & (1.512) & (1.630) & (1.586) & (1.505)   \\
\hline
\end{tabular}
}
\caption{{Coverage probabilities (CP) and coverage length (AL) for nominal $95\%$ confidence intervals.}}
\label{tab:ci}
\end{center}
\end{table}

\subsection{Illustrative example}
\label{sec:exm}

This example, primarily for illustration, uses the multivariate Nested-Error regression model (\ref{eqn:MNE}) and data from the posted land price data along the Keikyu train line from 1998 to 2001.
This train line connects the suburbs in the Kanagawa prefecture to the Tokyo metropolitan area.
Those who live in the suburbs in the Kanagawa prefecture take this line to work or study in Tokyo everyday.
Thus, it is expected that the land price depends on the distance from Tokyo.
The posted land price data are available for $53$ stations on the Keikyu train line, and we consider each station as a small area, namely, $m = 53$.

\medskip
For the $i$-th station, data of $n_i$ land spots are available, where $n_i$ varies around 4 and some areas have
only one observation.
For $i = 1,\ldots,m$, observations $\y_{ij}=(y_{ij1}, y_{ij2},y_{ij3})^\top$ denotes the difference between the value of the posted land price (Yen/1,000) for the unit meter squares of the $j$-th spot from 1998 to 2001, where $y_{ij1}$ is the a difference between 1998 and 1999, $y_{ij2}$ is the a difference between 1999 and 2000 and $y_{ij3}$ is the a difference between 2000 and 2001. 
As auxiliary variables, we use the data $(T_i, D_{ij}, FAR_{ij})$.
$T_i$ is the time to take from the nearby station $i$ to the Tokyo station around 8:30 in the morning, $D_{ij}$ is the value of geographical distance from the spot $j$ to the station $i$ and $FAR_{ij}$ denotes the floor-area ratio, or ratio of building volume to lot area of the spot $j$.
Then the regressor in the model (\ref{eqn:MNE}) is 
$$
\X_{ij}=
\left(
\begin{array}{cccccccccccc}
1 & FAR_{ij} & T_i & D_{ij} & 0 & 0 & 0 & 0 & 0 & 0 & 0 & 0 \\ 
0 & 0 & 0 & 0 &1 & FAR_{ij} & T_i & D_{ij} & 0 & 0 & 0 & 0 \\
0 & 0 & 0 & 0 & 0 & 0 & 0 & 0 &1 & FAR_{ij} & T_i & D_{ij}
\end{array}
\right)^\top.
$$

\medskip
The estimates of the covariance matrix $\bPsi$ and $\bSi$ are
$$
\bPsih=\begin{pmatrix}
43.4 & 27.3 & 28.4\\
27.3 & 33.4 & 20.4\\
28.4 & 20.4 & 28.5
\end{pmatrix}
\quad {\rm and}\quad
\bSih=\begin{pmatrix}
169.2 & 127.3 & 101.0\\
127.3 & 113.5 & 85.3\\
101.0 & 85.3 & 77.6
\end{pmatrix}.
$$
Thus, the estimated correlation coefficient of random effects $\bro=(\rho_{12},\rho_{13},\rho_{23})^\top$ is $(0.72,0.81,0.66)$, where $\rho_{ab}$ is the correlation coefficient of $v_a$ and $v_b$.
The estimates of the regression coefficients are $\bbeh=(-4.28,16.67,-1.79,-0.13,6.32,13.16,-2.34,-0.02,-4.22,11.16,-0.33,-0.06)^\top$.

\begin{figure}[htbp]
  \begin{center}
          \includegraphics[clip, width=9.9cm]{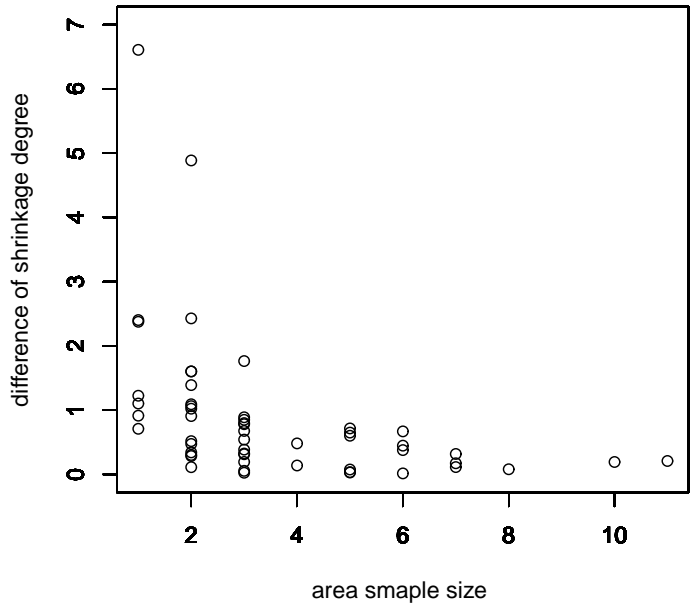}
          \hspace{1.6cm} 
    \caption{Plots of the difference of shrinkage degree against area sample size in MNER and NER.}
    \label{fig:sd}
  \end{center}
\end{figure}

All the estimated values of regression coefficients of $T_i$ and $D_{ij}$ are negative values which leads to the natural result that the $T_i$ and $D_{ij}$ have negative influence on $y_{ij}$, whose magnitude are almost unchanged for three years.
On the other hand, the magnitude of the influence of $FAR_{ij}$ on $y_{ij}$ decreases during the same time.
The obtained values of EBLUP for a difference between the posted land price data in  2000 and 2001 given in (\ref{eqn:EBLUP}) are give in Table \ref{tab:eb} for selected 15 areas.
To see the difference of predicted values of MNER and NER, Figure \ref{fig:sd} reports the difference between the degree of shrinkage, which is caluculated by $|{\rm dif}(\bthh^{EB})-{\rm dif}(\bthh^{uEB})|$ where ${\rm dif}(\bth)=|\yb_{ij3}-\bth_{ij3}|$.
It can be seen that the difference gets smaller as an area sample size $n_i$ gets larger.
This is because the smaple mean is reliable when $n_i$ is large, so that the sample mean does not be shrunk and the the degree of shrinkage of MNER and NER have almost no difference. 
In Table \ref{tab:eb}, we also provide the estimats of squared root of MSE (SMSE) given in (\ref{eqn:msemest}).
It is revealed from Table \ref{tab:eb} that SMSE of MNER is smaller than that of NER when $n_i$ is small.
On the other hand, SMSE of MNER is larger than that of NER when $n_i$ is large, particularly larger than $5$.
This is because the low accuracy in estimation of the covariance matrix $\bPsi$ and $\bSi$ has more adverse influence on prediction than the benefit from incorporating the small correlation into the estimation.
Table \ref{tab:ci} reports lower bounds (LB) and upper bounds (UB) of the $95\%$ confidence interval estimator of the difference between the value of the posted land price from 1998 to 2001, that is $\bel^\top\bth_a$ where $\bel=(1,1,1)^\top$, for selected 15 areas.

\begin{table}[htbp]
\begin{center}
\resizebox{11cm}{!} {
\begin{tabular}{ccccccc}
 \hline
& & sample & \multicolumn{2}{c}{MNER} & \multicolumn{2}{c}{NER} \\
area & $n_i$ & mean & EBLUP & SMSE &EBLUP & SMSE  \\
\hline
16 & 1& 3.0 & 2.94 & 4.93 & 9.66 & 5.03\\
17 & 1 & 23.0 & 29.70 & 4.85 & 27.32 & 4.96\\ 
31 & 1 & 75.0 & 38.93 & 4.74 & 41.33 & 4.85\\
21& 2 & 19.5 & 19.12 & 4.29 & 24.76 & 4.35\\ 
22 & 2 & 9.0 & 11.37 & 4.29 & 13.80 & 4.35\\
13 & 3 & 14.0 & 16.35 & 3.87 & 15.46 & 3.89\\ 
34 & 3 & 37.66 & 36.83 & 3.87 & 35.06 & 3.89\\
12 & 4 & 27.75 & 28.50 & 3.58 & 28.02 & 3.57\\ 
35 & 4 & 25.0 & 25.81 & 3.56 & 24.05 & 3.56\\
9 & 5 & 9.2 & 9.05 & 3.38 & 9.38 & 3.35\\
7 & 6 & 21.33 & 18.65 & 3.15 & 19.32 & 3.11\\
26 & 7 & 17.14 & 18.34 & 2.95 & 18.02 & 2.92\\
40 & 8 & 13.63 & 11.48 & 2.81 & 11.40 & 2.77\\   
52 & 10 & 6.0 & 6.33 & 2.58 & 6.52 & 2.54 \\
41 & 11 & 17.0 & 14.50 & 2.47 & 14.71 & 2.43 \\
\hline
\end{tabular}
}
\caption{{The estimated results for PLP Data for selected 15 areas.}}
\label{tab:eb}
\end{center}
\end{table}

\begin{table}[htbp]
\begin{center}
\resizebox{11cm}{!} {
\begin{tabular}{cccccc}
 \hline
area & $n_i$ & sample mean & EBLUP & LB & UB   \\
\hline
16 & 1& 73.0 & 48.99 & 14.97 & 83.02 \\
17 & 1 & 74.0 & 107.75 & 74.24 & 141.24 \\ 
31 & 1 & 261.0 & 136.79 & 103.93 & 169.65 \\
21& 2 & 132.5 & 109.78 & 80.78 & 138.78 \\ 
22 & 2 & 64.0 & 63.46 & 34.48 & 92.44 \\
13 & 3 & 50.0 & 58.86 & 32.61 & 85.11 \\ 
34 & 3 & 113.66 & 120.87 & 94.63 & 147.11 \\
12 & 4 & 103.5 & 104.60 & 80.20 & 129.00 \\ 
35 & 4 & 71.75 & 80.66 & 56.43 & 104.89 \\
9 & 5 & 23.2 & 21.91 & -1.28 & 45.11 \\
7 & 6 & 67.0 & 56.34 & 34.68 & 77.99 \\
26 & 7 & 65.86 & 69.62 & 49.25 & 89.99 \\
40 & 8 & 35.13 & 28.40 & 8.99 & 47.82 \\   
52 & 10 & 14.9 & 15.38 & -2.55 & 33.31 \\
41 & 11 & 47.55 & 38.91 & 21.75 & 56.06 \\
\hline
\end{tabular}
}
\caption{{$95\%$ confidence interval estimator for PLP Data for selected 15 areas.}}
\label{tab:ci}
\end{center}
\end{table}

\section{Proofs}
\label{sec:proof}

In this section, we use the notations $\bLa_i=\bPsi+n_i^{-1}\bSi$ and $\bLah_i=\bPsih+n_i^{-1}\bSih$ for $i=1, \ldots, m$.

\medskip
{\bf Proof of Theorem \ref{thm:consistency}}.\ \ 
We first show the results (1) and (2) for the estimators $\bPsih_1$ and $\bSih$.
Clearly, $E[\bSih]=\bSi$. 
To show that $E[\bPsih_1-\bPsi]=O(m^{-3/2})$, we begin by writing $\bPsih_0-\bPsi$ as
\begin{align*}
\bPsih_0-\bPsi=&
{1\over n}\sum_{i=1}^m\sum_{j=1}^{n_i}\{(\y_{ij}-\X_{ij}^\top\bbe)(\y_{ij}-\X_{ij}^\top\bbe)^\top - (\bPsi+\bSi)\}\\
&
+ {1\over n}\sum_{i=1}^m\sum_{j=1}^{n_i} \X_{ij}^\top(\bbeh^{OLS}-\bbe)(\bbeh^{OLS}-\bbe)^\top\X_{ij}
-{1\over n}\sum_{i=1}^m\sum_{j=1}^{n_i} (\y_{ij}-\X_{ij}^\top\bbe)(\bbeh^{OLS}-\bbe)^\top\X_{ij}
\\
&
-{1\over n}\sum_{i=1}^m\sum_{j=1}^{n_i} \X_{ij}^\top(\bbeh^{OLS}-\bbe)(\y_{ij}-\X_{ij}\bbe)
-(\bSih-\bSi),
\end{align*}
which yields the bias given in (\ref{eqn:PBias}).
Since $\bSih$ is unbiased, the bias ${\rm Bias}_{\bPsih_0}(\bPsi,\bSi)$ of $\bPsih_0$ is of order $O(m^{-1})$.
Using the results that $\bPsih_0-\bPsi=O_p(m^{-1/2})$ and $\bSih-\bSi=O_p(m^{-1/2})$, which will be shown below, we can see that $\bPsih_1=\bPsih_0-{\rm Bias}_{\bPsih_0}(\bPsih_0,\bSih)$ is a second-oder unbiased estimator of $\bPsi$.

\medskip
For (2), it is noted that $\bSih-\bSi$ is approximated as
\begin{align}
\bSih-\bSi=
{1\over N-m}\sum_{i=1}^m\Big\{\sum_{j=1}^{n_i} (\bep_{ij}-\bar \bep_i)(\bep_{ij}-\bar \bep_i)^\top-(n_i-1)\bSi\Big\}
+O_p(m^{-1}).
\label{eqn:Siap}
\end{align}
It is here noted that $\{\sum_{j=1}^{n_i} (\bep_{ij}-\bar \bep_i)(\bep_{ij}-\bar \bep_i)^\top-(n_i-1)\bSi\}/(N-m)$ for $i=1,\ldots,m$ are mutually  independent and $E\{\sum_{j=1}^{n_i} (\bep_{ij}-\bar \bep_i)(\bep_{ij}-\bar \bep_i)^\top-(n_i-1)\bSi\}/(N-m)=0$ for $i=1,\ldots,m$.
Then we can show that $\sqrt{m}(\bSih-\bSi)$ converges to a multivariate normal distributionthe because of the  finiteness of moments of normal random variables, which implies that $\bSih-\bSi=O_p(m^{-1/2})$.

Concerning $\bPsih_1-\bPsi=O_p(m^{-1/2})$, from the fact that ${\rm Bias}_{\bPsih_0}(\bPsih_0, \bSih)=O_p(m^{-1})$, it is sufficient to show that $\bPsih_0-\bPsi=O_p(m^{-1/2})$.
Then $\bPsih_0-\bPsi$ is approximated as
\begin{equation}
\bPsih_0-\bPsi=
{1\over N}\sum_{i=1}^m\Big\{\sum_{j=1}^{n_i}(\v_i+\bep_{ij})(\v_i+\bep_{ij})^\top - n_i(\bPsi+\bSi)\Big\}
-(\bSih_0-\bSi)
+O_p(m^{-1}).
\label{eqn:Psiap}
\end{equation}
It is here noted that $\{\sum_{j=1}^{n_i}(\v_i+\bep_{ij})(\v_i+\bep_{ij})^\top - n_i(\bPsi+\bSi)\}/N$ for $i=1,\ldots,m$ are mutually  independent and $E\{\sum_{j=1}^{n_i}(\v_i+\bep_{ij})(\v_i+\bep_{ij})^\top - n_i(\bPsi+\bSi)\}/N=0$ for $i=1,\ldots,m$.
Then we can show that $\sqrt{m}(\bPsih_0-\bPsi)$ converges to a multivariate normal distribution because of the  finiteness of moments of normal random variables, which implies that $\bPsih_0-\bPsi=O_p(m^{-1/2})$.

\medskip
We next prove (3) from the fact that $\sqrt{m}(\bPsih_1-\bPsi)=O_p(1)$.
The difference between $\bPsih$ and $\bPsih_1$ is in the case that $\bPsih_1$ is not nonegative definite.
Thus, we evaluare the probability $P(\a^\top\bPsih_1\a <0)$ for some $\a\in \Re^k$.
It is noted that the event $\a^\top\bPsih_1\a <0$ is equivalent to $-\sqrt{m}\a^\top(\bPsih_1-\bPsi)\a > \a^\top\bPsi\a$.
Using the Markov inequality, we observe that for any $\de>0$,
\begin{align*}
P(\a^\top\bPsih_1\a<0)
=& P(-\sqrt{m}\a^\top(\bPsih_1-\bPsi)\a>\sqrt{m}\a^\top\bPsi\a)
\\
\leq & P(|\sqrt{m}\a^\top(\bPsih_1-\bPsi)\a| > \sqrt{m}\a^\top\bPsi\a)
\\
\leq & E\Big[ \Big( {|\sqrt{m}\a^\top(\bPsih_1-\bPsi)\a|\over \sqrt{m}\a^\top\bPsi\a}\Big)^{2\de}\Big]=O(m^{-\de}),
\end{align*}
which proves (3) of Theorem \ref{thm:consistency}.

\medskip
Using the result (3) of Theorem \ref{thm:consistency}, we can show that $E[\bPsih]-\bPsi=O(m^{-3/2})$ and $\bPsih-\bPsi=O_p(m^{-1/2})$.

\medskip
Finally we verify that $\bbeh(\bPsih,\bSih)-\bbe=O_p(m^{-1/2})$.
Note that $\bbeh(\bPsih,\bSih)-\bbe$ is decomposed as $\{\bbeh(\bPsih,\bSih)-\bbeh(\bPsi,\bSi)\}+\{\bbeh(\bPsi,\bSi)-\bbe\}$.
For $\bbeh(\bPsi,\bSi)-\bbe$, it is noted that
\begin{align}
\bbeh(\bPsi,\bSi)-\bbe=
(\X^\top\D^{-1}\X)^{-1}\X^\top\D^{-1}(\y-\X\bbe).\non
\end{align}
Then, ${\rm Cov}(\bbeh(\bPsi,\bSi)-\bbe)=(\X^\top\D^{-1}\X)^{-1}=O(1/m)$ and this implies $\bbeh(\bPsi,\bSi)-\bbe=O_{p}(m^{-1/2})$.
We next evaluate $\bbeh(\bPsih,\bSih)-\bbeh(\bPsi,\bSi)$ as
\begin{align}
\bbeh(&\bPsih,\bSih)-\bbeh(\bPsi,\bSi)\non\\
=&(\X^\top\widehat\D^{-1}\X)^{-1}\X^\top\widehat\D^{-1}\y-(\X^\top\D^{-1}\X)^{-1}\X^\top\D^{-1}\y
\non\\
=&
(\X^\top\widehat\D^{-1}\X)^{-1}\X^\top(\widehat\D^{-1}-\D^{-1})\y
+\{(\X^\top\widehat\D^{-1}\X)^{-1}-(\X^\top\D^{-1}\X)^{-1}\}\X^\top\D^{-1}\y
\non\\
=&
I_{1}+I_{2},
\label{eqn:pbe1}
\end{align}
where $\widehat \D$ is obtained by replacing $\bSi$ and $\bPsi$ in $\D$ with $\bSih_0$ and $\bPsih_0$ respectively.
First, $I_{1}$ is written as
\begin{align}
I_{1}
=
-(\X^\top\widehat\D^{-1}\X)^{-1}\X^\top\widehat\D^{-1}(\widehat\D-\D)\D^{-1}\y,
\label{eqn:pbe2}
\end{align}
which is of order $O_p(m^{-1/2})$, because $(\X^\top\widehat\D^{-1}\X)=O_{p}(m)$ and $\X^\top\widehat\D^{-1}(\widehat\D-\D)\D^{-1}\y=O_{p}(m^{1/2})$.
Next, $I_{2}$ is rewritten as
\begin{align}
I_{2}=&
-(\X^\top\widehat\D^{-1}\X)^{-1}\X^\top(\widehat \D^{-1}-\D^{-1})\X(\X^\top\D^{-1}\X)^{-1}\X^\top\D^{-1}\y
\non\\
=&
(\X^\top\widehat\D^{-1}\X)^{-1}\X^\top\widehat \D^{-1}(\widehat \D-\D)\D^{-1}\X(\X^\top\D^{-1}\X)^{-1}\X^\top\D^{-1}\y
\label{eqn:pbe3}
\end{align}
which is of order $O_p(m^{-1/2})$, because $\X^\top\widehat\D^{-1}\X=O_{p}(m)$, $\X^\top\widehat \D^{-1}(\widehat \D-\D)\D^{-1}\X=O_{p}(m^{1/2})$, $\X^\top\D^{-1}\y=O_{p}(m)$. 
Thus, we have $\bbeh(\bPsih,\bSih)-\bbeh(\bPsi,\bSi)=O_p(m^{-1/2})$, and it is concluded that $\bbeh(\bPsih,\bSih)-\bbe=O_{p}(m^{-1/2})$.
\hfill$\Box$

\bigskip
{\bf Proof of Lemma \ref{lem:1}}.\ \ 
The covariance of $\y-\X\bbeh^{OLS}$ and $\bbeh(\bPsi,\bSi)$ is
\begin{align*}
E[&(\y-\X\bbeh^{OLS})(\bbeh(\bPsi,\bSi)-\bbe)^\top] \X^\top\D^{-1}\X\\
=&
E[ \{(\y-\X\bbe)-\X(\bbeh^{OLS}-\bbe)\}(\y-\X\bbe)^\top]\D^{-1}\X
\\
=&
(\I-\X(\X^\top\X)^{-1}\X^\top)E[ (\y-\X\bbe)(\y-\X\bbe)^\top]\D^{-1}\X
\\
=&
(\I-\X(\X^\top\X)^{-1}\X^\top)\D\D^{-1}\X=\zero,
\end{align*}
which implies that $\bbeh(\bPsi,\bSi)$ is independent of $\y-\X\bbeh^{OLS}$.
We next note that $\ybt=\Q\y$ and $\Xbt=\Q\X$ where $\Q={\rm block\ diag}(\P_1\otimes \I_k, \ldots, \P_m\otimes \I_k)$ for $\P_i=\I_{n_i}-n_i^{-1}\J_{n_i}$.
Then,
\begin{align*}
E[&(\ybt-\Xbt\bbet)(\bbeh(\bPsi,\bSi)-\bbe)^\top] \X^\top\D^{-1}\X\\
=&
E[ \{(\ybt-\Xbt\bbe)-\Xbt(\Xbt^{\top}\Xbt)^{-}\Xbt^{\top}(\ybt-\Xbt\bbe)^\top\}(\y-\X\bbe)^\top]\D^{-1}\X
\\
=&
\{\Q-\Q\X(\X^\top\Q^\top\Q\X)^{-}\X^\top\Q^\top\Q\}E[ (\y-\X\bbe)(\y-\X\bbe)^\top]\D^{-1}\X
\\
=&
\Q\X-\Q\X(\X^\top\Q^\top\Q\X)^{-}\X^\top\Q^\top\Q\X,
\end{align*}
which is equal to zero from the property of the generalized inverse.
Thus, $\bbeh(\bPsi,\bSi)$ is independent of $\ybt-\Xbt\bbet$, so that $\bbeh(\bPsi,\bSi)$ is independent of $\bSih$ and $\bPsih$.

It is also noted that
\begin{align*}
\bthh_a^{EB}&-\bthh_a(\bPsi,\bSi)
\\
=&
\c_a^\top \{\hbbe(\bPsih,\bSih)-\hbbe(\bPsi,\bSi) \}+ \bPsih\bLah_a^{-1}\{\yb_a - \Xb_a^\top\hbbe(\bPsih,\bSih)\}- \bPsi\bLa_a^{-1}\{\yb_a - \Xb_a^\top\hbbe(\bPsi,\bSi)\}
\\
=&
\c_a^\top \{\hbbe(\bPsih,\bSih)-\hbbe(\bPsi,\bSi) \}+ (\bPsih\bLah_a^{-1}-\bPsi\bLa_a^{-1})(\yb_a - \Xb_a^\top\bbeh^{OLS})
\\
&-\bPsih\bLah_a^{-1}\Xb_a^\top\{\bbeh(\bPsih,\bSih)-\bbeh^{OLS}\}
+\bPsi\bLa_a^{-1}\Xb_a^\top\{\hbbe(\bPsi,\bSi)-\bbeh^{OLS}\},
\end{align*}
which is a function of $\y-\X\bbeh^{OLS}$ and $\bSih$, because $\hbbe(\bPsi,\bSi))-\bbeh^{OLS}=(\X^\top\D^{-1}\X)^{-1}\X^\top\D^{-1}(\y-\X\bbeh^{OLS})-(\X^\top\X)^{-1}\X^\top(\y-\X\bbeh^{OLS})$.
Hence, $\bthh_a^{EB}-\bthh_a(\bPsi,\bSi)$ is independent of $\bbeh(\bPsi,\bSi)$.
\hfill$\Box$

\bigskip
{\bf Proof of Theorem \ref{thm:MSE}}.\ \ 
We shall prove that $E[\{\bthh_a^{EB}-\bthh_a(\bPsi,\bSi)\}\{\bthh_a^{EB}-\bthh_a(\bPsi,\bSi)\}^\top]=\G_{3a}(\bPsi,\bSi)+O_p(m^{-3/2})$.
It is observed that
\begin{align*}
\bthh_a^{EB}&-\bthh_a(\bPsi,\bSi)\\
=& (\bPsih\bLah_a^{-1}-\bPsi\bLa_a^{-1})(\yb_a-\Xb_a^\top\bbe)
+(\c_a^\top-\bPsih\bLah_a^{-1}\Xb_a^\top)\{\bbeh(\bPsih,\bSih)-\bbe\}
\\
&
-(\c_a^\top-\bPsi\bLa_a^{-1}\Xb_a^\top)\{\bbeh(\bPsi,\bSi)-\bbe\}.
\end{align*}
We can see that
\begin{align*}
(\bPsih&\bLah_a^{-1}-\bPsi\bLa_a^{-1})(\yb_a-\Xb_a^\top\bbe)\\
=& \big\{\bPsih(\bLah_a^{-1}-\bLa_a^{-1})+(\bPsih-\bPsi)\bLa_a^{-1}\big\}(\yb_a-\Xb_a^\top\bbe)
\\
=& \big\{(\I_{k}-\bPsih\bLah_a^{-1})(\bPsih-\bPsi)-\bPsih\bLah_a^{-1}n_a^{-1}(\bSih-\bSi)\big\}\bLa_a^{-1}(\yb_a-\Xb_a^\top\bbe)
\\
=& \big\{n_a^{-1}\bSi\bLa_a^{-1}(\bPsih-\bPsi)-\bPsi\bLa_a^{-1}n_a^{-1}(\bSih-\bSi)\big\}\bLa_a^{-1}(\yb_a-\Xb_a^\top\bbe)
+O_p(m^{-1})
\end{align*}
and
\begin{align*}
(\c_a^\top&-\bPsih\bLah_a^{-1}\Xb_a^\top)\{\bbeh(\bPsih,\bSih)-\bbe\}
\\
=& (\c_a^\top-\bPsi\bLa_a^{-1}\Xb_a^\top)\{\bbeh(\bPsih,\bSih)-\bbe\}
-(\bPsih\bLah_a^{-1}-\bPsi\bLa_a^{-1})\Xb_a^\top\{\bbeh(\bPsih,\bSih)-\bbe\}
\\
=&(\c_a^\top-\bPsi\bLa_a^{-1}\Xb_a^\top)\{\bbeh(\bPsih,\bSih)-\bbe\}
\\
& +\big\{n_a^{-1}\bSih\bLah_a^{-1}(\bPsih-\bPsi)-\bPsih\bLah_a^{-1}n_a^{-1}(\bSih-\bSi)\big\}\bLa_a^{-1}\Xb_a^\top\{\bbeh(\bPsih,\bSih)-\bbe\}
\\
=& \big\{\c_a^\top-\bPsi\bLa_a^{-1}\Xb_a^\top\big\}\{\bbeh(\bPsih,\bSih)-\bbe\} +O_p(m^{-1}).
\end{align*}
Thus, we have
\begin{align*}
\bthh_a^{EB}-\bthh_a(\bPsi,\bSi)
=&
\big\{n_a^{-1}\bSi\bLa_a^{-1}(\bPsih-\bPsi)-\bPsi\bLa_a^{-1}n_a^{-1}(\bSih-\bSi)\big\}\bLa_a^{-1}(\yb_a-\Xb_a^\top\bbe)\\
&
+ (\c_a^\top-\bPsi\bLa_a^{-1}\Xb_a^\top)\{\bbeh(\bPsih,\bSih)-\bbe\} +O_p(m^{-1})
\\
=&I_1+I_2 + O_p(m^{-1}). \quad \text{(say)}
\end{align*}
For $I_2$, it is noted that
\begin{align*}
\bbeh(\bPsih,\bSih)-\bbeh(\bPsi,\bSi)
=& 
\big\{(\X^\top\widehat\D^{-1}\X)^{-1}-(\X^\top\D^{-1}\X)^{-1}\big\}\X^\top\widehat\D^{-1}(\y-\X\bbe)
\\
&+(\X^\top\D^{-1}\X)^{-1}\X^\top(\widehat\D^{-1}-\D^{-1})(\y-\X\bbe)
\\
=&I_{21}+I_{22},  \quad \text{(say)}.
\end{align*}
We can evaluate $I_{21}$ as
\begin{align*}
I_{21}=(\X^\top\D^{-1}\X)^{-1}\X^\top\widehat\D^{-1}(\widehat\D-\D)\D^{-1}\X\{\bbeh(\bPsih,\bSi)-\bbe\}= O_p(m^{-1}),
\end{align*}
because $\X^\top\D^{-1}\X=O(m)$, $\X^\top\widehat\D^{-1}(\widehat\D-\D)\D^{-1}\X =O_p(m^{1/2})$ and $\bbeh(\bPsih,\bSih)-\bbe=O_p(m^{-1/2})$ from Theorem \ref{thm:consistency} (2).
We next estimate $I_{22}$ as
\begin{align*}
I_{22}=&
- (\X^\top\D^{-1}\X)^{-1}
\Big\{\sum_{i=1}^m \X_i^\top\A(\bPsih,\bSih)\X_i\Big\}\\
&\times \Big\{\sum_{i=1}^m \X_i^\top \A(\bPsih,\bSih)\X_i\Big\}^{-1}\sum_{i=1}^m \X_i^\top\A(\bPsih,\bSih)(\y_i-\X_i\bbe),
\end{align*}
where
$$
\A(\bPsih,\bSih)=(\J_{n_i}\otimes \bPsih+\I_{n_i}\otimes \bSih)^{-1}\big\{\J_{n_i}\otimes (\bPsih-\bPsi)+\I_{n_i}\otimes (\bSih-\bSi)\big\}(\J_{n_i}\otimes \bPsi+\I_{n_i}\otimes \bSi)^{-1}.
$$
It can be seen that $I_{22}=O_p(m^{-1})$ from the same arguments as in $I_{21}$.
Thus, it follows that $I_2=O_p(m^{-1})$.

From equations (\ref{eqn:Siap}) and (\ref{eqn:Psiap}),
\begin{align}
\bthh_a^{EB}&-\bthh_a(\bPsi,\bSi)
\non\\
=& \Big\{n_a^{-1}\bSi\bLa_a^{-1}(\bPsih-\bPsi)-\bPsi\bLa_a^{-1}n_a^{-1}(\bSih-\bSi)\Big\}\bLa_a^{-1}(\yb_a-\Xb_a^\top\bbe)
+O_p(m^{-1})
\label{eqn:ebap} \\
=& \Big[{n_a^{-1} \over N}\bSi\bLa_a^{-1}\sum_{i=1}^m \sum_{j=1}^{n_i} \big\{ (\v_i+\bep_{ij})(\v_i+\bep_{ij})^\top-(\bPsi+\bSi)\big\}
\non\\
& -{n_a^{-1} \over N-m}(\bPsi+\bSi)\bLa_a^{-1}\sum_{i=1}^m\sum_{j=1}^{n_i}\big\{(\bep_{ij}-\bar\bep_i)(\bep_{ij}-\bar\bep_i)^\top-(1-n_i^{-1})\bSi  \big\}\Big]
\non\\
&\times \bLa_a^{-1} (\yb_a-\Xb_a^\top\bbe)+O_p(m^{-1})
\non\\
=& {n_a^{-1} \over N}\bSi\bLa_a^{-1}\sum_{i=1}^m n_i\big\{ (\v_i+\bar\bep_i)(\v_i+\bar\bep_i)^\top-\bLa_i\big\}\bLa_a^{-1} (\v_a+\bar\bep_a)
\non\\
&-{n_a^{-1} \over N(N-m)}(N\bPsi+m\bSi)\bLa_a^{-1}\sum_{i=1}^m \sum_{j=1}^{n_i} \big\{ (\bep_{ij}-\bar\bep_i)(\bep_{ij}-\bar\bep_i)^\top-(1-n_i^{-1})\bSi\big\}
\non\\
&\quad \times \bLa_a^{-1} (\v_a+\bar\bep_a)+O_p(m^{-1})
\non\\
=& \A_1(\y) - \A_2(\y)+O_p(m^{-1}), \quad \text{(say)}.
\label{eqn:ebap1}
\end{align}
Note that $E[\A_1(\y)\{\A_2(\y)\}^\top]=\zero$ because $\v_i+\bar\bep_i$ is independent of $\bep_{ij}-\bar\bep_i$. 
Hence, we need to evaluate $E[\A_1(\y)\{\A_1(\y)\}^\top]$ and $E[\A_2(\y)\{\A_2(\y)\}^\top]$. 

Concerning $E[\A_1(\y)\{\A_1(\y)\}^\top]$, it can be seen that
\begin{align*}
E\Big[&\sum_{i=1}^m n_i \big\{ (\v_i+\bar\bep_i)(\v_i+\bar\bep_i)^\top-\bLa_i\big\}\bLa_a^{-1} (\v_a+\bar\bep_a)(\v_a+\bar\bep_a)^\top\bLa_a^{-1}
\sum_{j=1}^m n_j \big\{ (\v_j+\bar\bep_j)(\v_j+\bar\bep_j)^\top-\bLa_j\big\}\Big]
\\
=&
\sum_{i=1}^m n_i^2 E\Big[\big\{ (\v_i+\bar\bep_i)(\v_i+\bar\bep_i)^\top-\bLa_i\big\}\bLa_a^{-1} (\v_a+\bar\bep_a)(\v_a+\bar\bep_a)^\top\bLa_a^{-1}
\big\{ (\v_i+\bar\bep_i)(\v_i+\bar\bep_i)^\top-\bLa_i\big\}\Big]
\\
=&
\sum_{i\neq a} n_i^2 E\Big[\big\{ (\v_i+\bar\bep_i)(\v_i+\bar\bep_i)^\top-\bLa_i\big\}\bLa_a^{-1} \big\{ (\v_i+\bar\bep_i)(\v_i+\bar\bep_i)^\top-\bLa_i\big\}\Big]
+O(1)
\\
=&
\sum_{i\neq a} n_i^2\big\{ \bLa_i\bLa_a^{-1}\bLa_i+\tr(\bLa_a^{-1}\bLa_i)\bLa_i\big\}
+O(1)
=
\sum_{i=1}^m n_i^2\big\{ \bLa_i\bLa_a^{-1}\bLa_i+\tr(\bLa_a^{-1}\bLa_i)\bLa_i\big\}
+O(1),
\end{align*}
so that we have
\begin{equation}
E[\A_1(\y)\{\A_1(\y)\}^\top]
={n_a^{-2} \over N^2}\bSi\bLa_a^{-1}\sum_{i=1}^m n_i^2\big\{ \bLa_i\bLa_a^{-1}\bLa_i
+\tr(\bLa_a^{-1}\bLa_i)\bLa_i\big\}\bLa_a^{-1}\bSi
+O(m^{-2}).
\label{eqn:ebap2}
\end{equation}

Concerning the evaluation of $E[\A_2(\y)\{\A_2(\y)\}^\top]$, let $\W=\sum_{i=1}^m \sum_{j=1}^{n_i}  (\bep_{ij}-\bar\bep_i)(\bep_{ij}-\bar\bep_i)^\top$ for simplicity.
Then, $\W$ has the Wishart distribution $\Wc_k(N-m, \bSi)$.
Because $\v_a+\bar\bep_a$ is independent of $\bep_{ij}-\bar\bep_i$, it follows that
\begin{align*}
E\Big[&\sum_{i=1}^m \sum_{j=1}^{n_i} \big\{ (\bep_{ij}-\bar\bep_i)(\bep_{ij}-\bar\bep_i)^\top-(1-n_i^{-1})\bSi\big\}\bLa_a^{-1} (\v_a+\bar\bep_a)(\v_a+\bar\bep_a)^\top\bLa_a^{-1} 
\\
&\times \sum_{k=1}^m \sum_{\ell=1}^{n_k} \big\{ (\bep_{k\ell}-\bar\bep_k)(\bep_{k\ell}-\bar\bep_k)^\top-(1-n_k^{-1})\bSi\big\}\Big]
\\
=&
E\big[ \{\W-(N-m)\bSi\}\bLa_a^{-1} (\v_a+\bar\bep_a)(\v_a+\bar\bep_a)^\top\bLa_a^{-1} 
\{\W-(N-m)\bSi\}\big]
\\
=&
E\big[ \{\W-(N-m)\bSi\}\bLa_a^{-1}\{\W-(N-m)\bSi\}\big].
\end{align*}
From the properties of the Wishart distribution, it is noted that $E[\W]=(N-m)\bSi$ and $E[\W\bLa_a^{-1}\W]=(N-m)(N-m+1)\bSi\bLa_a^{-1}\bSi+(N-m)\tr(\bLa_a^{-1}\bSi)\bSi$.
Thus,
\begin{align}
E[\A_2(\y)\{\A_2(\y)\}^\top]
=& {n_a^{-2} \over N^2(N-m)^2}(N\bPsi+m\bSi)\bLa_a^{-1}E\big[ \{\W-(N-m)\bSi\}\bLa_a^{-1}\{\W-(N-m)\bSi\}\big]
\non\\
&\times \bLa_a^{-1}(N\bPsi+m\bSi)
\non\\
=& {n_a^{-2} \over N^2(N-m)}(N\bPsi+m\bSi)\bLa_a^{-1}\Big\{\bSi\bLa_a^{-1}\bSi+\tr(\bLa_a^{-1}\bSi)\bSi\Big\}\bLa_a^{-1}(N\bPsi+m\bSi).
\label{eqn:ebap3}
\end{align}
Combining (\ref{eqn:ebap1}), (\ref{eqn:ebap2}) and (\ref{eqn:ebap3}) gives the expression in (\ref{eqn:G3}).
\hfill$\Box$

\bigskip
{\bf Proof of Theorem \ref{thm:msemest}}.\ \ From (2) in Theorem \ref{thm:consistency}, it is sufficient to show this approximation for $\bPsih$.
We can rewrite $G_{1a}(\bPsih)$ as 
\begin{align*}
\G_{1a}(\bPsih,\bSih)
=& n_a^{-1}\bPsih(\bPsih+n_a^{-1}\bSih)^{-1}\bSih
\\
=&
\G_{1a}(\bPsi,\bSi) +  n_a^{-2}\bSi\bLa_a^{-1}(\bPsih-\bPsi)\bLa_a^{-1}\bSi
+n_a^{-1}\bPsi\bLa_a^{-1}(\bSih-\bSi)\bLa_a^{-1}\bPsi\\
&-n_a^{-2}\bSi\bLa_a^{-1}(\bPsih-\bPsi)\bLa_a^{-1}(\bPsih-\bPsi)\bLa_a^{-1}\bSi
-n_a^{-2}\bPsi\bLa_a^{-1}(\bSih-\bSi)\bLa_a^{-1}(\bSih-\bSi)\bLa_a^{-1}\bPsi\\
&+n_a^{-2}\bSi\bLa_a^{-1}(\bPsih-\bPsi)\bLa_a^{-1}(\bSih-\bSi)\bLa_a^{-1}\bPsi
+n_a^{-2}\bPsi\bLa_a^{-1}(\bSih-\bSi)\bLa_a^{-1}(\bPsih-\bPsi)\bLa_a^{-1}\bSi
\\
&+O_p(m^{-3/2}),
\end{align*}
which implies that
\begin{align*}
\G_{1a}&(\bPsi,\bSi) - E[\G_{1a}(\bPsih,\bSih)]\\
=
&E\Big[n_a^{-2}\bSi\bLa_a^{-1}(\bPsih-\bPsi)\bLa_a^{-1}(\bPsih-\bPsi)\bLa_a^{-1}\bSi
+n_a^{-2}\bPsi\bLa_a^{-1}(\bSih-\bSi)\bLa_a^{-1}(\bSih-\bSi)\bLa_a^{-1}\bPsi\\
&-n_a^{-2}\bSi\bLa_a^{-1}(\bPsih-\bPsi)\bLa_a^{-1}(\bSih-\bSi)\bLa_a^{-1}\bPsi
-n_a^{-2}\bPsi\bLa_a^{-1}(\bSih-\bSi)\bLa_a^{-1}(\bPsih-\bPsi)\bLa_a^{-1}\bSi\Big]
\\
&+O(m^{-3/2}),
\end{align*}
because $\bPsih$ is second-order unbiased and $\bSih$ is unbiased.
On the other hand, from (\ref{eqn:ebap}), it follows that
\begin{align*}
E[&\{\bthh_a^{EB}-\bthh_a(\bPsi,\bSi)\}\{\bthh_a^{EB}-\bthh_a(\bPsi,\bSi)\}^\top]\\
=&E\Big[\big\{n_a^{-1}\bSi\bLa_a^{-1}(\bPsih-\bPsi)-\bPsi\bLa_a^{-1}n_a^{-1}(\bSih-\bSi)\big\}\bLa_a^{-1}(\bar \y_a-\bar \X_a^\top\bbe)
\\
&\times (\bar \y_a-\bar \X_a^\top\bbe)^\top \bLa_a^{-1}\big\{(\bPsih-\bPsi)\bLa_a^{-1}n_a^{-1}\bSi-n_a^{-1}(\bSih-\bSi)\bLa_a^{-1}\bPsi\big\}\Big]
+O(m^{-3/2})
\\
=&E\Big[\big\{n_a^{-1}\bSi\bLa_a^{-1}(\bPsih-\bPsi)-\bPsi\bLa_a^{-1}n_a^{-1}(\bSih-\bSi)\big\}\bLa_a^{-1}\big\{(\bPsih-\bPsi)\bLa_a^{-1}n_a^{-1}\bSi-n_a^{-1}(\bSih-\bSi)\bLa_a^{-1}\bPsi\big\}\Big]\\
&+O(m^{-3/2}).
\end{align*}
Thus, we have
$$
\G_{1a}(\bPsi,\bSi) - E[\G_{1a}(\bPsih,\bSih)]
=
E[\{\bthh_a^{EB}-\bthh_a(\bPsi,\bSi)\}\{\bthh_a^{EB}-\bthh_a(\bPsi,\bSi)\}^\top]
+O(m^{-3/2}),
$$
which yields $\G_{1a}(\bPsi,\bSi) - E[\G_{1a}(\bPsih,\bSih)]=\G_{3a}(\bPsi,\bSi)+O(m^{-3/2})$.
Since $\G_{3a}(\bPsi,\bSi)=O(m^{-1})$, one gets 
$$
\G_{1a}(\bPsi,\bSi) = E[\G_{1a}(\bPsih,\bSih) + \G_{3a}(\bPsih,\bSih)]+O(m^{-3/2}),
$$
and Theorem \ref{thm:msemest} is established.
\hfill$\Box$

\bigskip
{\bf Proof of Theorem \ref{thm:ci}}.\ \ 
The proof is done along the line given in Diao et al. (2014).
Let $\P_\X=\I_k-\X(\X^\top\X)^{-1}\X^\top$.
From Lemma \ref{lem:1}, $\bel^\top(\bthh_a^{EB}-\bthh_a(\bPsi, \bSi))$ is a function of $\P_\X\y$ and is independent of $\bel^\top(\bthh_a(\bPsi, \bSi)-\bth_a)$ given $\P_\X\y$.
It is noted that $E[\bel^\top(\bthh_a(\bPsi, \bSi)-\bth_a)(\bthh_a(\bPsi, \bSi)-\bth_a)^\top\bel]=\bel^\top E[(\bthh_a(\bPsi, \bSi)-\btht_a(\bbe, \bPsi,\bSi))(\bthh_a(\bPsi, \bSi)-\btht_a(\bbe, \bPsi,\bSi))^\top+(\btht_a(\bbe, \bPsi,\bSi)-\bth_a)(\btht_a(\bbe, \bPsi,\bSi)-\bth_a)^\top]\bel=\bel^\top(\G_{1a}(\bPsi, \bSi)+\G_{2a}(\bPsi, \bSi))\bel$, where $\G_{1a}(\bPsi, \bSi)$ and $\G_{2a}(\bPsi, \bSi)$ are given in (\ref{eqn:G12}).
Then the conditional distribution of $\bel^\top(\bthh_a^{EB}-\bth_a)$ given $\P_\X \y$ is
\begin{align}
\bel^\top(\bthh_a^{EB}-\bth_a) | \P_\X \y \sim \Nc_k (\bel^\top(\bthh_a^{EB}-\bthh_a(\bPsi, \bSi)), \bel^\top\H_a(\bPsi, \bSi)\bel),
\label{eqn:cd}
\end{align}
where $\H_a(\bPsi, \bSi)=\G_{1a}(\bPsi, \bSi)+\G_{2a}(\bPsi, \bSi)$.
This implies that
\begin{align*}
P\Big(&{\bel^\top(\bthh_a^{EB}-\bth_a) \over \{\bel^\top{\rm msem}(\bthh_a^{EB})\bel\}^{1/2}} \leq z\Big)
\\
=&
E\Big[P\Big({\bel^\top(\bthh_a^{EB}-\bth_a)-\bel^\top(\bthh_a^{EB}-\bthh_a(\bPsi, \bSi)) \over \{\bel^\top\H_a(\bPsi, \bSi)\bel\}^{1/2}}\\
&\qquad\quad\leq {\{\bel^\top{\rm msem}(\bthh_a^{EB})\bel\}^{1/2}z-\bel^\top(\bthh_a^{EB}-\bthh_a(\bPsi, \bSi))  \over \{\bel^\top\H_a(\bPsi, \bSi)\bel\}^{1/2} } | \P_\X\y\Big)\Big]
\\
=&
E\Big[\Phi\Big( {\{\bel^\top{\rm msem}(\bthh_a^{EB})\bel\}^{1/2}z-\bel^\top(\bthh_a^{EB}-\bthh_a(\bPsi, \bSi))  \over \{\bel^\top\H_a(\bPsi, \bSi)\bel\}^{1/2} }  \Big)\Big].
\end{align*}
Thus, it is observed that 
\begin{align*}
P\Big(& \bel^\top\bthh_a^{EB}-\{\bel^\top{\rm msem}(\bthh_a^{EB})\bel\}^{1/2}z \leq \bel^\top\bth_a \leq \bel^\top\bthh_a^{EB}+\{\bel^\top{\rm msem}(\bthh_a^{EB})\bel\}^{1/2}z \Big)
\\
=&
E\Big[\Phi\Big( {\{\bel^\top{\rm msem}(\bthh_a^{EB})\bel\}^{1/2}z-\bel^\top(\bthh_a^{EB}-\bthh_a(\bPsi, \bSi))  \over \{\bel^\top\H_a(\bPsi, \bSi)\bel\}^{1/2} }  \Big)
\\
&\quad
-\Phi\Big( {-\{\bel^\top{\rm msem}(\bthh_a^{EB})\bel\}^{1/2}z-\bel^\top(\bthh_a^{EB}-\bthh_a(\bPsi, \bSi))  \over \{\bel^\top\H_a(\bPsi, \bSi)\bel\}^{1/2} }  \Big)   \Big]
\\
=&
E\Big[\Phi(r_{1a}-r_{2a})-\Phi(-r_{1a}-r_{2a})\Big]
=
E\Big[\Phi(r_{1a}+r_{2a})+\Phi(r_{1a}-r_{2a})\Big]-1,
\end{align*}
where 
\begin{align*}
r_{1a}=&\{\bel^\top{\rm msem}(\bthh_a^{EB})\bel\}^{1/2}z/\{\bel^\top\H_a(\bPsi, \bSi)\bel\}^{1/2},\\
r_{2a}=&\bel^\top(\bthh_a^{EB}-\bthh_a(\bPsi, \bSi))/\{\bel^\top\H_a(\bPsi, \bSi)\bel\}^{1/2}.
\end{align*}
By the Taylor series expansion, for $r_{1a}^* \in (r_{1a},r_{1a}+r_{2a})$ and $r_{1a}^{**} \in (r_{1a},r_{1a}-r_{2a})$, we have
\begin{align}
\Phi(r_{1a}+r_{2a})+\Phi(r_{1a}-r_{2a})=2\Phi(r_{1a})+r_{2a}^2\phi^{(1)}(r_{1a})+{1 \over 24}r_{2a}^4(\phi^{(3)}(r_{1a}^*)+\phi^{(3)}(r_{1a}^**)),
\label{eqn:ci1}
\end{align}
where $\phi^{(1)}(\cdot)$ and $\phi^{(3)}(\cdot)$ are the first and third derivatives of the standard normal density $\phi(\cdot)$.
The Taylor series expansion is also used to get
\begin{align}
\{\bel^\top{\rm msem}(\bthh_a^{EB})\bel\}^{1/2}
=&\{\bel^\top {\rm MSEM}(\bthh_a^{EB}) \bel+\bel^\top{\rm msem}(\bthh_a^{EB})\bel-\bel^\top {\rm MSEM}(\bthh_a^{EB}) \bel\}^{1/2}
\non\\
\begin{split}
=&
\{\bel^\top {\rm MSEM}(\bthh_a^{EB}) \bel\}^{1/2}\Big(1+{ \bel^\top{\rm msem}(\bthh_a^{EB})\bel-\bel^\top {\rm MSEM}(\bthh_a^{EB}) \bel \over 2\bel^\top {\rm MSEM}(\bthh_a^{EB}) \bel }
\\
& - { (\bel^\top{\rm msem}(\bthh_a^{EB})\bel-\bel^\top {\rm MSEM}(\bthh_a^{EB})\bel )^2 \over 8\{\bel^\top {\rm MSEM}(\bthh_a^{EB}) \bel\}^2 }+\cdots  \Big).
\end{split}
\label{eqn:ci2}
\end{align}
We evaluate the expectation of the first term $2\Phi(r_{1a})$ in (\ref{eqn:ci1}).
By the Taylor series expansion, for $z^* \in (z,r_{1a})$, it is seen that
\begin{align}
\Phi(r_{1a})-\Phi(z)=(r_{1a}-z)\phi(z)+{(r_{1a}-z)^2 \over 2}\phi^{(1)}(z)+{(r_{1a}-z)^3 \over 6}\phi^{(2)}(z)+{(r_{1a}-z)^4 \over 24}\phi^{(1)}(z^*).
\label{eqn:ci25}
\end{align}
From (\ref{eqn:ci2}), we can evaluate $E[r_{1a}]$ as
\begin{align}
E[r_{1a}]
=&
\Big({\bel^\top{\rm MSEM}(\bthh_a^{EB})\bel \over \bel^\top\H_a(\bPsi, \bSi)\bel}\Big)^{1/2}z\Big[ 1+ {E[ \bel^\top{\rm msem}(\bthh_a^{EB})\bel-\bel^\top {\rm MSEM}(\bthh_a^{EB}) \bel] \over 2\bel^\top {\rm MSEM}(\bthh_a^{EB}) \bel } 
\non\\
\begin{split}
&- { E[(\bel^\top{\rm msem}(\bthh_a^{EB})\bel-\bel^\top {\rm MSEM}(\bthh_a^{EB})\bel )^2] \over 8\{\bel^\top {\rm MSEM}(\bthh_a^{EB}) \bel\}^2 }+ { E[(\bel^\top{\rm msem}(\bthh_a^{EB})\bel-\bel^\top {\rm MSEM}(\bthh_a^{EB})\bel )^3] \over 16\{\bel^\top {\rm MSEM}(\bthh_a^{EB}) \bel\}^3 }
\\
&+{1\over 16}{5 \over 8}E\Big[ \int_{\bel^\top{\rm msem}(\bthh_a^{EB})\bel}^{\bel^\top {\rm MSEM}(\bthh_a^{EB})\bel} \{\bel^\top {\rm MSEM}(\bthh_a^{EB})\bel\}^{-1/2}x^{-7/2}(\bel^\top{\rm msem}(\bthh_a^{EB})\bel-\bel^\top {\rm MSEM}(\bthh_a^{EB})\bel )^3 \Big] \Big].
\end{split}
\label{eqn:ci3}
\end{align}
Since ${\rm MSEM}(\bthh_a^{EB})=O(1)$ and ${\rm msem}(\bthh_a^{EB})$ is a second order unbised estimator of ${\rm MSEM}(\bthh_a^{EB})$, the second term in the bracket of (\ref{eqn:ci3}) is of order $o(m^{-1})$.
From Lemma \ref{lem:ci}, the moments of higher than three are of order $o(m^{-1})$, and we have $E[\{\bel^\top{\rm msem}(\bthh_a^{EB})\bel-\bel^\top {\rm MSEM}(\bthh_a^{EB})\bel \}^2]=V(\bthh_a^{EB})+o(m^{-1})$.
Then,  using 
$$\Big({\bel^\top{\rm MSEM}(\bthh_a^{EB})\bel \over \bel^\top\H_a(\bPsi, \bSi)\bel}\Big)^{1/2}=1+{1\over 2}{\bel^\top\G_{3a}(\bthh_a^{EB})\bel \over \bel^\top\H_a(\bPsi, \bSi)\bel}+o(m^{-1}),$$
we have
\begin{align}
E[r_{1a}]-z=&\Big({\bel^\top{\rm MSEM}(\bthh_a^{EB})\bel \over \bel^\top\H_a(\bPsi, \bSi)\bel}\Big)^{1/2}z\Big[ 1- { V(\bthh_a^{EB})\over 8\{\bel^\top {\rm MSEM}(\bthh_a^{EB}) \bel\}^2 }\Big]-z+o(m^{-1})
\non\\
=&
\Big\{ \Big({\bel^\top{\rm MSEM}(\bthh_a^{EB})\bel \over \bel^\top\H_a(\bPsi, \bSi)\bel}\Big)^{1/2}-1- \Big({\bel^\top{\rm MSEM}(\bthh_a^{EB})\bel \over \bel^\top\H_a(\bPsi, \bSi)\bel}\Big)^{1/2}{ V(\bthh_a^{EB}) \over 8\{\bel^\top {\rm MSEM}(\bthh_a^{EB}) \bel\}^2 }\Big\}z
\non\\
&+o(m^{-1})
\non\\
=&
\Big\{{1\over 2}{\bel^\top\G_{3a}(\bthh_a^{EB})\bel \over \bel^\top\H_a(\bPsi, \bSi)\bel}-{V(\bthh_a^{EB}) \over 8\{\bel^\top {\rm MSEM}(\bthh_a^{EB}) \bel\}^2 }\Big\}z+o(m^{-1}),
\label{eqn:ci4}
\end{align}
because $\G_{3a}(\bthh_a^{EB})$ and $V(\bthh_a^{EB})$ are of order $O(m^{-1})$.

Since $E[r_{1a}^2]=E[\bel^\top{\rm msem}(\bthh_a^{EB})\bel z^2/\bel^\top\H_a(\bPsi, \bSi)\bel=z^2\bel^\top{\rm MSEM}(\bthh_a^{EB})\bel/\bel^\top\H_a(\bPsi, \bSi)\bel+o(m^{-1})$ and $E[(r_{1a}-z)^2]=E[r_{1a}^2]-2zE[r_{1a}-z]-z^2$, it is observed that
\begin{align}
&E[(r_{1a}-z)^2]
\non\\
=&
z^2{\bel^\top{\rm MSEM}(\bthh_a^{EB})\bel \over \bel^\top\H_a(\bPsi, \bSi)\bel}-2z\Big(\Big({\bel^\top{\rm MSEM}(\bthh_a^{EB})\bel \over \bel^\top\H_a(\bPsi, \bSi)\bel}\Big)^{1/2}z\Big[ 1- { E[(\bel^\top{\rm msem}(\bthh_a^{EB})\bel-\bel^\top {\rm MSEM}(\bthh_a^{EB})\bel )^2] \over 8\{\bel^\top {\rm MSEM}(\bthh_a^{EB}) \bel\}^2 }\Big]
\non\\
&-z\Big)-z^2+o(m^{-1})
\non\\
=&
\Big\{ \Big({\bel^\top{\rm MSEM}(\bthh_a^{EB})\bel \over \bel^\top\H_a(\bPsi, \bSi)\bel}\Big)^{1/2}-1\Big)^2 + \Big({\bel^\top{\rm MSEM}(\bthh_a^{EB})\bel \over \bel^\top\H_a(\bPsi, \bSi)\bel}\Big)^{1/2} { V(\bthh_a^{EB})\over 4\{\bel^\top {\rm MSEM}(\bthh_a^{EB}) \bel\}^2 } \Big\}z^2+o(m^{-1})
\non\\
=&
{ V(\bthh_a^{EB})\over 4\{\bel^\top {\rm MSEM}(\bthh_a^{EB}) \bel\}^2 } z^2+o(m^{-1}).
\label{eqn:ci5}
\end{align}
This implies that $r_{1a}-z=O_p(m^{-1/2})$.
Then, the expectation of the third and forth terms of (\ref{eqn:ci25}) is of order $O_(m^{-3/2})$.

\medskip
We evaluate the expectation of the second term in (\ref{eqn:ci1}).
Since $E[r_{2a}^2]=\bel^\top\G_{3a}\bel/\bel^\top\H_a\bel$ and $E[r_{1a}]-z$ are of order $O(m^{-1})$, we have
\begin{align}
E[r_{2a}^2\phi^{(1)}(r_{1a})]=E[r_{2a}^2\phi^{(1)}(z)]+o(m^{-1})={\bel^\top\G_{3a}\bel \over \bel^\top\H_a\bel}\phi^{(1)}(z)+o(m^{-1}).
\label{eqn:ci6}
\end{align}

\medskip
Since $r_{2a}^4=O(m^{-2})$ and $E[r_{1a}^*]-z=O(1)$, the expectation of the third term in (\ref{eqn:ci1}) is of order $O(m^{-2})$.

\medskip
Combining (\ref{eqn:ci1}), (\ref{eqn:ci25}), (\ref{eqn:ci4}),  (\ref{eqn:ci5}) and  (\ref{eqn:ci6}) gives
\begin{align*}
P\Big(& \bel^\top\bthh_a^{EB}-\{\bel^\top{\rm msem}(\bthh_a^{EB})\bel\}^{1/2}z \leq \bel^\top\bth_a \leq \bel^\top\bthh_a^{EB}+\{\bel^\top{\rm msem}(\bthh_a^{EB})\bel\}^{1/2}z \Big)
\\
=&
2\Phi(z)+2\Big\{ {1\over 2}{\bel^\top\G_{3a}(\bthh_a^{EB})\bel \over \bel^\top\H_a(\bPsi, \bSi)\bel}-{ V(\bthh_a^{EB}) \over 8\{\bel^\top {\rm MSEM}(\bthh_a^{EB}) \bel\}^2 }\Big\}z\phi(z)
\\
&+ { V(\bthh_a^{EB}) \over 4\{\bel^\top {\rm MSEM}(\bthh_a^{EB}) \bel\}^2 } z^2\phi^{(1)}(z)
+{\bel^\top\G_{3a}\bel \over \bel^\top\H_a\bel}\phi^{(1)}(z)-1+o(m^{-1})
\\
=&
2\Phi(z)-1- { V(\bthh_a^{EB})\over 4\{\bel^\top {\rm MSEM}(\bthh_a^{EB}) \bel\}^2 } (z^3+z)\phi(z)+o(m^{-1}),
\end{align*}
which establishes Theorem \ref{thm:ci}.
\hfill$\Box$

\bigskip
{\bf Proof of Lemma \ref{lem:ci}}.\ \ 
It is noted that
\begin{align*}
\G_{1a}&(\bPsih,\bSih)-\G_{1a}(\bPsi,\bSi)
=
n_a^{-2}\bSi\bLa_a^{-1}(\bPsih-\bPsi)\bLa_a^{-1}\bSi
+n_a^{-1}\bPsi\bLa_a^{-1}(\bSih-\bSi)\bLa_a^{-1}\bPsi
+O_p(m^{-1}),
\end{align*}
and $\G_{2a}(\bPsih, \bSih)$, $\G_{2a}(\bPsi, \bSi)$, $\G_{3a}(\bPsi, \bSi)$ and $\G_{3a}(\bPsih, \bSih)$ are of order $O(m^{-1})$.
Thus, 
\begin{align*}
E[&(\bel^\top{\rm msem}(\bthh_a^{EB})\bel-\bel^\top {\rm MSEM}(\bthh_a^{EB})\bel )^2]
\\
=&
\bel^\top E(n_a^{-2}\bSi\bLa_a^{-1}(\bPsih-\bPsi)\bLa_a^{-1}\bSi
+n_a^{-1}\bPsi\bLa_a^{-1}(\bSih-\bSi)\bLa_a^{-1}\bPsi)\bel\bel^\top
\\
&\times(n_a^{-2}\bSi\bLa_a^{-1}(\bPsih-\bPsi)\bLa_a^{-1}\bSi
+n_a^{-1}\bPsi\bLa_a^{-1}(\bSih-\bSi)\bLa_a^{-1}\bPsi)\bel+o(m^{-1}),
\end{align*}
and the moment of $\bel^\top{\rm msem}(\bthh_a^{EB})\bel-\bel^\top {\rm MSEM}(\bthh_a^{EB})\bel$ of order higher than three is of order $o(m^{-1})$.

From equations (\ref{eqn:Siap}) and (\ref{eqn:Psiap}), it follows that
\begin{align*}
n_a^{-2}\bSi&\bLa_a^{-1}(\bPsih-\bPsi)\bLa_a^{-1}\bSi
+n_a^{-1}\bPsi\bLa_a^{-1}(\bSih-\bSi)\bLa_a^{-1}\bPsi)
\\
=&
{n_a^{-2} \over N}\bSi\bLa_a^{-1}\sum_{i=1}^m n_i\big\{ (\v_i+\bar\bep_i)(\v_i+\bar\bep_i)^\top-\bLa_i\big\}\bLa_a^{-1} \bSi
\non\\
&- \Big[{n_a^{-2} m\over N(N-m)}\bSi\bLa_a^{-1}\sum_{i=1}^m \sum_{j=1}^{n_i} \big\{ (\bep_{ij}-\bar\bep_i)(\bep_{ij}-\bar\bep_i)^\top-(1-n_i^{-1})\bSi\big\}\bLa_a^{-1}\bSi
\non\\
&\qquad -{n_a^{-1} \over N-m}\bPsi\bLa_a^{-1}\sum_{i=1}^m \sum_{j=1}^{n_i} \big\{ (\bep_{ij}-\bar\bep_i)(\bep_{ij}-\bar\bep_i)^\top-(1-n_i^{-1})\bSi\big\}\bLa_a^{-1}\bPsi\Big]+O_p(m^{-1})
\non\\
=& \B_1(\y) - \B_2(\y)+O_p(m^{-1}), \quad \text{(say)}.
\end{align*}
Since $\v_i+\bar\bep_i$ is independent of $\bep_{ij}-\bar\bep_i$, it is seen that $E[\B_1(\y)\bel\bel^\top\{\B_2(\y)\}^\top]=\zero$. 
Thus, we shall evaluate $E[\B_1(\y)\bel\bel^\top\{\B_1(\y)\}^\top]$ and $E[\B_2(\y)\bel\bel^\top\{\B_2(\y)\}^\top]$. 
For the proofs, we can use the same arguments as in (\ref{eqn:ebap1}), (\ref{eqn:ebap2}) and (\ref{eqn:ebap3}).

Concerning $E[\B_1(\y)\bel\bel^\top\{\B_1(\y)\}^\top]$, it is observed that
\begin{align*}
E\Big[&\sum_{i=1}^m n_i\big\{ (\v_i+\bar\bep_i)(\v_i+\bar\bep_i)^\top-\bLa_i\big\}\bLa_a^{-1} \bSi\bel\bel^\top\bSi\bLa_a^{-1}\sum_{j=1}^m n_j\big\{ (\v_j+\bar\bep_j)(\v_j+\bar\bep_j)^\top-\bLa_j\big\} \Big]
\\
=&
\sum_{i=1}^m n_i^2 E\Big[\big\{ (\v_i+\bar\bep_i)(\v_i+\bar\bep_i)^\top-\bLa_i\big\}\bLa_a^{-1} \bSi\bel\bel^\top\bSi\bLa_a^{-1}\big\{ (\v_i+\bar\bep_i)(\v_i+\bar\bep_i)^\top-\bLa_i\big\} \Big]
\\
=&
\sum_{i=1}^m n_i^2 \{ \bLa_i\bLa_a^{-1} \bSi\bel\bel^\top\bSi\bLa_a^{-1}\bLa_i+ \tr(\bLa_a^{-1} \bSi\bel\bel^\top\bSi\bLa_a^{-1}\bLa_i)\bLa_i\},
\end{align*}
so that we have
\begin{align}
E[\B_1(\y)\bel\bel^\top\{\B_1(\y)\}^\top]={n_a^{-4} \over N^2}\sum_{i=1}^m n_i^2 \{ \bLa_i\bLa_a^{-1} \bSi\bel\bel^\top\bSi\bLa_a^{-1}\bLa_i+ \tr(\bLa_a^{-1} \bSi\bel\bel^\top\bSi\bLa_a^{-1}\bLa_i)\bLa_i\}.
\label{eqn:lci1}
\end{align}

Concerning $E[\B_2(\y)\bel\bel^\top\{\B_2(\y)\}^\top]$, recall that $\W=\sum_{i=1}^m \sum_{j=1}^{n_i}  (\bep_{ij}-\bar\bep_i)(\bep_{ij}-\bar\bep_i)^\top$ has the Wishart distribution $\Wc_k(N-m, \bSi)$.
Then, we have
\begin{align}
E[&\B_2(\y)\bel\bel^\top\{\B_2(\y)\}^\top]
\non\\
=&
E\Big\{{n_a^{-2} m\over N(N-m)}\bSi\bLa_a^{-1} \big\{ \W-(N-m)\bSi\big\}\bLa_a^{-1}\bSi
-{n_a^{-1} \over N-m}\bPsi\bLa_a^{-1}\big\{ \W-(N-m)\bSi\big\}\bLa_a^{-1}\bPsi\Big\}
\non\\
&\times \bel\bel^\top \Big\{{n_a^{-2} m\over N(N-m)}\bSi\bLa_a^{-1} \big\{ \W-(N-m)\bSi\big\}\bLa_a^{-1}\bSi
-{n_a^{-1} \over N-m}\bPsi\bLa_a^{-1}\big\{ \W-(N-m)\bSi\big\}\bLa_a^{-1}\bPsi\Big\}
\non\\
\begin{split}
=&
{n_a^{-4} m^2\over N^2(N-m)}\bSi\bLa_a^{-1}\Big\{ \bSi\bLa_a^{-1} \bSi\bel\bel^\top\bSi\bLa_a^{-1}\bSi+\tr (\bLa_a^{-1} \bSi\bel\bel^\top\bSi\bLa_a^{-1}\bSi)\bSi  \Big\}\bLa_a^{-1}\bSi
\\
&+{n_a^{-2} \over N-m}\bPsi\bLa_a^{-1}\Big\{ \bSi\bLa_a^{-1} \bPsi\bel\bel^\top\bPsi\bLa_a^{-1}\bSi+\tr (\bLa_a^{-1} \bPsi\bel\bel^\top\bPsi\bLa_a^{-1}\bSi)\bSi  \Big\}\bLa_a^{-1}\bPsi
\\
&-{n_a^{-3} m\over N(N-m)}\bSi\bLa_a^{-1}\Big\{ \bSi\bLa_a^{-1} \bSi\bel\bel^\top\bPsi\bLa_a^{-1}\bSi+\tr (\bLa_a^{-1} \bSi\bel\bel^\top\bPsi\bLa_a^{-1}\bSi)\bSi \Big\}\bLa_a^{-1}\bPsi
\\
&-{n_a^{-3} m\over N(N-m)}\bPsi\bLa_a^{-1}\Big\{ \bSi\bLa_a^{-1} \bPsi\bel\bel^\top\bSi\bLa_a^{-1}\bSi+\tr (\bLa_a^{-1} \bPsi\bel\bel^\top\bSi\bLa_a^{-1}\bSi)\bSi \Big\}\bLa_a^{-1}\bSi.
\end{split}
\label{eqn:lci2}
\end{align}
Multipling $\bel$ by (\ref{eqn:lci1}) and (\ref{eqn:lci2}) from both sides, we get the expression given in (\ref{eqn:vmse}).
\hfill$\Box$

\section*{Acknowledgments}
Research of the second author was supported in part by Grant-in-Aid for Scientific Research  (15H01943 and 26330036) from Japan Society for the Promotion of Science.

\end{document}